\documentclass{article}%
\usepackage{amsmath}
\usepackage{amsfonts}
\usepackage{amssymb}
\usepackage{graphicx}
\usepackage{color}
\usepackage{hyperref}%
\providecommand{\U}[1]{\protect\rule{.1in}{.1in}}
\newtheorem{theorem}{Theorem}[section]

\newtheorem{corollary}[theorem]{Corollary}

\newtheorem{definition}[theorem]{Definition}

\newtheorem{lemma}[theorem]{Lemma}

\newtheorem{question}[theorem]{Question}
\newtheorem{proposition}[theorem]{Proposition}
\newtheorem{remark}[theorem]{Remark}

\newenvironment{proof}[1][Proof]{\noindent\textbf{#1.} }{\ \rule{0.5em}{0.5em}}
\begin{document}

\date{}
\title{Compactness of products of block Hankel and Toeplitz operators}
\author{Caixing Gu, Meng Li and Pan Ma}
\maketitle

\begin{abstract}
Motivated by the Sarason problem on the products of Hankel and Toeplitz
operators on analytic function spaces, we characterize the compactness of
products of block Hankel and Toeplitz operators on the vector-valued Hardy
space of the unit disk via harmonic extension of the symbols and Douglas
algebras generated by the symbols. Additionally, we provide a complete answer
to the question of when the product of a block Hankel operator and a block
Toeplitz operator is a block Hankel operator.

\end{abstract}

\section{Introduction}

Let $L^{2}$ be the space of square-integrable functions on the unit circle
$\mathbb{T}$ with respect to the normalized Lebesgue measure. Let $H^{2}$ be
the Hardy space on the open unit disk $\mathbb{D}$. $L^{\infty}$ and
$H^{\infty}$ are the algebras of bounded functions in $L^{2}$ and $H^{2}$
respectively. Let $B(E\mathcal{)}$ be the algebra of bounded linear operators
on a complex separable Hilbert space $E$. $L_{E}^{2}$ and $H_{E}^{2}$ denote
$E$-valued $L^{2}$ and $H^{2}$ spaces. $L_{B(E)}^{\infty}$ and $H_{B(E)}%
^{\infty}$ are operator-valued $L^{\infty}$ and $H^{\infty}$ algebras. Let
$M_{n\times n}$ denote the space of $n\times n$ complex matrices. $L_{n\times
n}^{\infty}$ and $H_{n\times n}^{\infty}$ denote spaces of size $n\times n$
with entries being $L^{\infty}$ and $H^{\infty}$ functions, respectively. The
theory of operator-valued analytic functions can be found in \cite[Section V]{NFBK}.

Let $\Phi\in L_{B(E)}^{\infty},$ then the block Toeplitz operator $T_{\Phi}$
from $H_{E}^{2}$ into $H_{E}^{2}$ is defined by%
\[
T_{\Phi}h=P\left[  \Phi h\right]  ,h\in H_{E}^{2},
\]
where $P$ denotes the projection from $L_{E}^{2}$ to $H_{E}^{2}.$ Let $T_{z}%
$($=T_{zI_{E}}$) denote the shift operator on $H_{E}^{2}$ for some $E$ which
the context will make it clear. We will also use $S_{E}$ and $S_{E}^{\ast}$ to
denote this shift operator $T_{z}$ and the backward shift $T_{z}^{\ast}$ on
$H_{E}^{2}.$ In the scalar Hardy space $H^{2},$ we just use $S$ and $S^{\ast}%
$. The block Toeplitz operator $T_{\Phi}=A$ is characterized by the operator
equation $S_{E}^{\ast}AS_{E}=A.$ It is well known that the block Toeplitz
operators $T_{\Phi}$ are bounded and are not compact except $\Phi=0$ for
$\Phi\in L_{B(E)}^{\infty}$ in \cite[Page 103-104]{Peller}. Block Toeplitz operators and
(block) Toeplitz determinants (i.e., determinants of sections of (block)
Toeplitz operators) are not only important in a wide range of fields such as
the famous Ising model \cite{DIK1, DIK2}, but also represent one of the most
significant classes of non-selfadjoint operators. They provide a fascinating
example of the fruitful interaction between operator theory, the theory of
Banach algebras and function theory; see \cite{BS2, Douglas, Zhu2007}.
H. Widom \cite{Widom} characterized the asymptotic behavior of block Toeplitz
determinants with smooth symbols.
Inspired by applications to entanglement entropy, E. Basor, T. Ehrhardt and
J.Virtanen \cite{Basor} studied the asymptotic behavior of block Toeplitz
determinants based on a new method of computing the Fredholm index of Toeplitz
operators with piecewise continuous matrix-valued symbols.

For $f\in L_{E}^{2}$ and $\Phi\in L_{B(E)}^{\infty},$ let%
\[
\widetilde{f}(z)=f(\overline{z}),\widetilde{\Phi}(z)=\Phi(\overline{z}).
\]
Let $J$ be defined on $L_{E}^{2}$ by
\[
Jf(z)=\overline{z}f(\overline{z})=\overline{z}\widetilde{f},f(z)\in L_{E}%
^{2}.
\]
Then $J$ maps $\overline{zH_{E}^{2}}$ onto $H_{E}^{2}$ and $J$ maps $H_{E}%
^{2}$ onto $\overline{zH_{E}^{2}}.$ Furthermore, $J$ is a unitary operator,
\[
J^{\ast}=J,J^{2}=I,\text{ }JQ=PJ \text{ and }JP=QJ,
\]
where $Q:=I-P$ denotes the projection from $L_{E}^{2}$ to $\overline
{zH_{E}^{2}}:=L_{E}^{2}\ominus H_{E}^{2}.$ The block Hankel operator $H_{\Phi
}$ from $H_{E}^{2}$ into $H_{E}^{2}$ is defined by%
\[
H_{\Phi}h=JQ\left[  \Phi h\right]  =PJ\left[  \Phi h\right]  ,h(z)\in
H_{E}^{2}.
\]
The block Hankel operator $H_{\Phi}=A$ is characterized by the operator
equation $AS_{E}=S_{E}^{\ast}A$. It is easy to see that $H_{\Phi}$ is bounded
for $\Phi\in L_{B(E)}^{\infty}.$ The block Hankel operator $H_{\Phi}$ is
compact on $H_{E}^{2}$ if and only if $\Phi\in H_{B(E)}^{\infty}+C(B(E))$,
where $C(B(E))$ denotes the $B(E)$-valued continuous function on the unit
circle $\mathbb{T}$ in \cite{Hartman, Page}. Block Hankel operators are not
only used in control, approximation and prediction theory, but also are
important in non-commutative analysis, see \cite{Par, Peller}. For example, G.
Pisier \cite{pis} solved the famous problem of whether a polynomially bounded
operator on Hilbert space must be similar to a contraction with the help of
block Hankel operators.

In this paper $E$ will be finite dimensional, so we will use $E$ and
$\mathbb{C}^{n}$ interchangeably. In this case, if $\{e_{i},1\leq i\leq n\}$
is any orthonormal basis of $E,$ then%
\begin{equation}
C_{E}:=I-S_{E}S_{E}^{\ast}=\sum_{i=1}^{n}e_{i}\otimes e_{i}, \label{defce}%
\end{equation}
where $e_{i}\otimes e_{i}$ denotes the rank one operator. That is, $C_{E}$ is
the projection from $H_{E}^{2}$ to the constant functions $E.$ Often
$\{e_{i},1\leq i\leq n\}$ is the standard basis of $\mathbb{C}^{n}.$ We note
that sometimes $T_{\Phi}$ and $H_{\Phi}$ could be unbounded operators, where
$\Phi\in L_{B(E)}^{2}$; that is, $\Phi$ is a matrix-valued function whose
entries belong to $L^{2}.$ In this case, an operator equation involving
$T_{\Phi}$ and $H_{\Phi}$ is valid when it acts on polynomials or $H^{\infty}$
functions; see for example \cite{GuH}, where the product of unbounded operator
$T_{\Phi}H_{\Psi}$ can be interpreted by using bilinear form $\left\langle
T_{\Phi}H_{\Psi}h,g\right\rangle :=\left\langle H_{\Psi}h,T_{\Phi^{\ast}%
}g\right\rangle $ for $h,g$ being polynomials or $H^{\infty}$ functions.

In \cite{Sarason problem}, D. Sarason proposed the problem of characterizing the pairs of functions $f,g\in H^2$ such that the operator $T_{f}T_{\overline{g}}$ is bounded on $H^2$ and presented a necessary condition \eqref{Sarason conjecture} which leads to the following conjecture referred to as \textbf{Sarason's conjecture}:

{\sl For two functions $f,g\in H^2,$ $T_{f}T_{\overline{g}}$ is bounded if and only if 
\begin{equation}\label{Sarason conjecture}
\sup_{z\in \mathbb D}\widehat{|f|^2}(z)\widehat{|g|^2}(z)<\infty,
\end{equation}
where $\widehat{h}$ is defined as the Poisson extension of $h\in L^1(\mathbb T).$}

D. Zheng \cite{Zheng} proved that (\ref{Sarason conjecture}) with $2$ replaced by $2+\epsilon$ is sufficient by the distribution function inequality. K. Stroethoff and D. Zheng \cite{strzheng2} showed that for $ f, g \in H^2 $, $ T_f T_{\overline{g}} $ is bounded and invertible on $ H^2 $ if and only if \eqref{Sarason conjecture} holds and $ \inf\{ |f(z)||g(z)| : z \in \mathbb{D} \} > 0 $.
In addition, using a number of product identities for Toeplitz operators, K. Stroethoff and D. Zheng \cite{strzheng2} proved that $ T_f T_{\overline{g}} $ is bounded and Fredholm on $ H^2 $ if and only if \eqref{Sarason conjecture}
holds and the function $|f||g|$  is
 bounded away from zero on $\mathbb{T}$.

From the perspective of operator theory, Sarason’s Toeplitz product problem and the
associated Sarason conjecture make perfect sense for the Bergman spaces and the Fock spaces. Natural companion to Sarason’s original Toeplitz product problem is to characterize the boundedness and compactness of product of Toeplitz and Hankel operators, product of Hankel operators on analytic function spaces which was called \textbf{Sarason’s Ha-plitz Product Problem} by K. Zhu \cite{Zhu23}.
K. Stroethoff and D. Zheng made a series of outstanding and sustained contributions to characterize the Sarason’s Ha-plitz product problem on Bergman spaces \cite{strzheng1,strzheng2,strzheng3,strzheng4,strzheng5}. While the Sarason’s conjecture are false on the Hardy space \cite{Nazarov} and the Bergman space \cite{Alpore}, it turns out that the Fock space version of Sarason’s conjecture is true \cite{chpazh}. P. Ma, F. Yan, D. Zheng and K. Zhu \cite{mayanzhengzhu1} solved the natural conjectures on product of Hankel operators on Fock space and completely characterized the boundedness and compactness of product of Toeplitz and Hankel operators on Fock space in \cite{mayanzhengzhu2}. There has been additional work in recent years about analogs of Sarason’s problems and conjectures for Toeplitz operators on weighted Fock spaces \cite{BYZ,CWXC,FLL}. More results can be found in the recent surveys \cite{Virtanen22,Zhu23}.

Due to the non-commutative property of the product of
matrices, the characterization of the product of block Hankel operators and
block Toeplitz operators on the vector-valued Hardy space is more difficult
than in the scalar case \cite{GuZheng}. Some new phenomena arise in the
product of block Hankel operators on the vector-valued Hardy spaces of the
unit disk. For instance, it is well known that the product of two scalar
Hankel operators on the Hardy space is zero if and only if one of the
operators must be zero \cite{Brown Halmos}. However, this property does not
hold for block Hankel operators with matrix-valued symbols; see
\cite{Curto2023}. The well-known Brown-Halmos theorem does not hold for the
product of block Toeplitz operators \cite{GuZheng}.

Motivated by the Sarason problem \cite{Sarason problem} and the application to
entanglement entropy \cite{Basor}, the primary focus of this paper is
exploring products of block Hankel and Toeplitz operators on the vector-valued
Hardy space of the unit disk $\mathbb{D}$. In this paper, we will present a
complete solution to the following natural question.

\begin{question}
Assume $E=\mathbb{C}^{n}$ and $\Phi,\Psi\in L_{B(E)}^{\infty}.$ When is the
product of block Hankel operators and block Toeplitz operators $H_{\Phi
}T_{\Psi}$ compact on the vector-valued Hardy space $H_{E}^{2}$ of the unit
disk $\mathbb{D}?$
\end{question}

In order to state the first main result, we need to recall the harmonic
extension. Let $g(w)\in L^{1}=L^{1}(\mathbb{T)}$, and $g(z)$ denote its
harmonic extension at $z\in\mathbb{D}$, via Poisson integral%
\[
g(z)=\int_{\mathbb{T}}g(w)\frac{1-\left\vert z\right\vert ^{2}}{\left\vert
1-\overline{z}w\right\vert ^{2}}d\sigma(w),
\]
where $d\sigma(w)$ is the normalized Lebesgue measure on $\mathbb{T}$. In
fact, we often identify a function $g(w)\in L^{2}=L^{2}(\mathbb{T)}$ with its
harmonic extension $g(z).$ In this sense, $H^{2}$ is the subspace of harmonic
functions which happen to be analytic. If $g\in L_{E}^{1}$ or $g$ is a
matrix-valued function, then the harmonic extension of $g$ is defined
entrywise. If we write $\Phi\in L_{B(E)}^{\infty}$ as $\Phi=\Phi_{+}+\Phi_{-}$
with $\Phi_{+}=P(\Phi)$ in $H_{E}^{2}$ and $\Phi_{-}=Q(\Phi)$ in $\overline{zH_{E}^{2}}$, then
the answer to the above question can be expressed by the harmonic extensions
of the symbols $\Phi$ and $\Psi$.

\begin{theorem}
\label{thm harmonic extension} Assume $E=\mathbb{C}^{n}$ and $\Phi,\Psi\in
L_{B(E)}^{\infty}.$ Then $H_{\Phi}T_{\Psi}$ is compact on the vector-valued
Hardy space $H_{E}^{2}$ if and only if%
\[
\lim_{\left\vert z\right\vert \rightarrow1}trace\left[  \left\vert \Phi
_{-}-\Phi_{-}(z)\right\vert ^{2}(z)\cdot\left\vert \widetilde{\Psi_{-}}^{\ast
}-\widetilde{\Psi_{-}}^{\ast}(\overline{z})\right\vert ^{2}(\overline
{z})\right]  =0
\]
and
\[
\lim_{\left\vert z\right\vert \rightarrow1}trace\left[  \left\vert \left(
\Phi_{-}\Psi_{+}\right)  _{-}-\left(  \Phi_{-}\Psi_{+}\right)  _{-}(z)+\left(
\left[  \Phi_{-}-\Phi_{-}(z)\right]  \Psi_{-}(z)\right)  \right\vert
^{2}(z)\right]  =0.
\]

\end{theorem}

In the scalar case, we get a new result that complements Theorem
\ref{ChengChu15} below.

\begin{corollary}
Let $\varphi,\psi\in L^{\infty}.$ Then $H_{\varphi}T_{\psi}$ is compact on the
Hardy space $H^{2}$ if and only if%
\[
\lim_{\left\vert z\right\vert \rightarrow1}\left\vert \varphi_{-}-\varphi
_{-}(z)\right\vert ^{2}(z)\left(  \left\vert \psi_{-}-\psi_{-}(z)\right\vert
^{2}(z)\right)  =0
\]
and
\[
\lim_{\left\vert z\right\vert \rightarrow1}\left[  \left\vert \left(  \phi
_{-}\psi_{+}\right)  _{-}-\left(  \phi_{-}\psi_{+}\right)  _{-}(z)+\left(
\left[  \phi_{-}-\phi_{-}(z)\right]  \psi_{-}(z)\right)  \right\vert ^{2}(z)
\right]  =0.
\]

\end{corollary}

Next we need to recall the definition of a Douglas algebra. A closed
subalgebra between $H^{\infty}$ and $L^{\infty}$ is called a Douglas algebra
\cite{Chang, Marshall}. $H^{\infty}[\varphi]$ denotes the Douglas algebra of
$L^{\infty}$ generated by $H^{\infty}$ and $\varphi$. The second main result
of our paper is as follows.

\begin{theorem}
\label{maina}Assume $E=\mathbb{C}^{n}$ and $\Phi,\Psi\in L_{B(E)}^{\infty}.$
Then $H_{\Phi}T_{\Psi}$ is compact on the vector-valued Hardy space $H_{E}%
^{2}$ if and only if
\[
\bigcap_{A\in\left(  M_{n}\right)  _{2^{2n}}}H^{\infty}\left[  \Phi
(I-A),A\Psi,\Phi A\Psi\right]  \subset H_{n\times n}^{\infty}+C_{n\times n},
\]
where \[
\bigcap_{A \in (M_{n})_{2^{2n}}} 
   H^{\infty}\!\bigl[\Phi(I-A),\, A\Psi,\, \Phi A\Psi\bigr]
\]denotes the intersection, over all $A \in (M_{n})_{2^{2n}}$, 
of the Douglas algebras generated by $H^{\infty}$ together with $\Phi(I-A)$, $A\Psi$ and $\Phi A\Psi$.
\end{theorem}

There are three motivations of our results. The first motivation of our
results is the beautiful theorem of Axler, Chang and Sarason in \cite{ACS} and
Volberg in \cite{Volberg}, which concerns the product of scalar Hankel
operators on the Hardy space $H^{2}.$ In addition, D. Xia and D. Zheng  \cite{ZhengXia1, ZhengXia2}  completely characterized compact products of three Hankel operators on $H^2$ by the Hoffman's abstract Hardy theory on the support set in \cite{Hoffmanacta}. C. Gu \cite{GUIEOT} gave a necessary and sufficient condition for the product of several Hankel operators on $H^2$ to be compact.

\begin{theorem}
\label{ACS,Volberg} \cite{ACS,Volberg}Let $\varphi,\psi\in L^{\infty}.$ Then
$H_{\varphi}^{\ast}H_{\psi}$ is compact on the Hardy space $H^{2}$ if and only
if
\[
H^{\infty}[{\varphi}]\cap H^{\infty}[\psi]\subset H^{\infty}+C.
\]

\end{theorem}

The second motivation of our results is the following theorem given by Zheng
\cite{Zheng}, characterizing the compactness of the product of two Hankel
operators, complementing Theorem \ref{ACS,Volberg} above.

\begin{theorem}
\cite{Zheng} Let $\varphi,\psi\in L^{\infty}.$ Then $H_{\varphi}^{\ast}%
H_{\psi}$ is compact on the Hardy space $H^{2}$ if and only if
\[
\lim_{\left\vert z\right\vert \rightarrow1}\left\vert \varphi_{-}-\varphi
_{-}(z)\right\vert ^{2}(z)\cdot\left\vert \psi_{-}-\psi_{-}(z)\right\vert
^{2}(z)=0.
\]
Here $\left\vert \varphi_{-}-\varphi_{-}(z)\right\vert ^{2}(z)$ is the
harmonic extension at $z$ of $\left\vert \varphi_{-}(w)-\varphi_{-}%
(z)\right\vert ^{2}$ as a function of $w$ in $L^{1}(\mathbb{T)}$.
\end{theorem}

The third motivation of our results is Chu's theorem \cite{ChengChu15} which
gave the necessary and sufficient conditions for the compactness of the
product of a scalar Hankel operator and a scalar Toeplitz operator on the
Hardy space.

\begin{theorem}
\label{ChengChu15} \cite{ChengChu15} Let $\varphi,\psi\in L^{\infty}.$ Then
$H_{\varphi}T_{\psi}$ is compact on the Hardy space $H^{2}$ if and only if%
\[
H^{\infty}[\varphi]\cap H^{\infty}[\psi,\varphi\psi]\subset H^{\infty}+C.
\]

\end{theorem}

Now, we outline the paper. In Section \ref{Section 2}, we first introduce some
facts about Hankel and Toeplitz operators, which will be used later. In order
to study the compactness of the product of block Hankel operators and block
Toeplitz operators, it is necessary to study the question of when their
product is zero or another block Hankel operator. We answer the question of
when $H_{\Phi}T_{\Psi}$ is a block Hankel operator in Theorem \ref{hankel},
and provide a more explicit description of the symbols $\Phi$ and $\Psi$ in
Theorem \ref{Theorem Huw}. In Section \ref{Section 3}, we present the
necessary and sufficient conditions for the compactness of a sum of odd
products of block Hankel operators and block Toeplitz operators, as stated in
Theorem \ref{main}. In addition, Theorem \ref{Theorem 3.8} describes the
compactness of the product of block Hankel operators and block Toeplitz
operators through the harmonic extension of their symbols. In Section
\ref{Section 4}, to get Theorem \ref{maina}, we will prove several equivalent
conditions for the compactness of $H_{\Phi}T_{\Psi}$ in Theorem
\ref{main theorem TFAE}. The proof of Theorem \ref{main theorem TFAE} is
complicated and long. The proof applies the theory of maximal ideal space of
$H^{\infty}$, Douglas algebras, localization techniques as pioneered by Axler,
Chang, Sarason, Volberg, Zheng and others, it also use a number of powerful
algebraic insights on manipulating products of block Hankel and Toeplitz
operators. In Section \ref{Section 5}, we will present some applications of
Theorem \ref{main theorem TFAE} with concrete examples to highlight the
complexity and challenges involved in Theorem \ref{maina}.

\section{Product of block Hankel and Toeplitz operators}

\label{Section 2}

The study of the compactness of the product of a Hankel operator and a
Toeplitz operator naturally leads to the question of when their product is
zero or another Hankel operator. It turns out this question is considered by
the first author in \cite{GuJOT} where a partial answer is given, see
\cite{Ding}. Later, Lee \cite{Lee} proved under what conditions the finite sum
of products of Hankel and Toeplitz operators equals zero. C. Chu \cite{ChengChu20} gave the necessary and sufficient conditions for the compactness of the
finite sum of products of Hankel and Toeplitz operators. The tools we will
develop to study the compactness of the product of block Hankel and Toeplitz
operators render an answer to this question. We first collect a few facts
about block Hankel operators and block Toeplitz operators with brief proofs.

\begin{lemma}
Let $\Phi,\Psi\in L_{B(E)}^{\infty}.$ Then $T_{\Phi}^{\ast}=T_{\Phi^{\ast}},$
$H_{\Phi}^{\ast}=H_{\widetilde{\Phi}^{\ast}},$ and%
\begin{align}
T_{\Phi\Psi}  &  =T_{\Phi}T_{\Psi}+H_{\widetilde{\Phi}}H_{\Psi},\label{t1}\\
H_{\Phi\Psi}  &  =H_{\Phi}T_{\Psi}+T_{\widetilde{\Phi}}H_{\Psi}. \label{h1}%
\end{align}

\end{lemma}

\begin{proof}
With respect to the decomposition $L_{E}^{2}=H_{E}^{2}\oplus\overline
{zH_{E}^{2}},$ the multiplication operator $M_{\Phi}$ can be written as%
\[
M_{\Phi}=\left[
\begin{array}
[c]{cc}%
T_{\Phi} & H_{\widetilde{\Phi}}J\\
JH_{\Phi} & JT_{\widetilde{\Phi}}J
\end{array}
\right]  .
\]
Then writing $M_{\Phi}M_{\Psi}=M_{\Phi}{}_{\Psi}$ as a $2\times2$ matrix
operator equation will prove the lemma.
\end{proof}

\begin{corollary}
\label{cort2h2}Let $\Phi\in L_{B(E)}^{\infty}$ and $\Psi\in H_{B(E)}^{\infty
}.$ Then%
\begin{align}
T_{\Phi}T_{\Psi}  &  =T_{\Phi\Psi},\label{t2}\\
H_{\Phi}T_{\Psi}  &  =H_{\Phi\Psi}. \label{h2}%
\end{align}
Similarly, $T_{\Psi}^{\ast}T_{\Phi}=T_{\Psi^{\ast}\Phi}$ and $T_{\widetilde
{\Psi}}H_{\Phi}=H_{\Psi\Phi}.$
\end{corollary}

\begin{proof}
By (\ref{t1}), $T_{\Phi\Psi}=T_{\Phi}T_{\Psi}+H_{\widetilde{\Phi}}H_{\Psi}.$
Note $H_{\Psi}=0$ since $\Psi\in H_{B(E)}^{\infty}.$ So $T_{\Phi}T_{\Psi
}=T_{\Phi\Psi}.$
\end{proof}

The formula (\ref{t2}) is related to Brown-Halmos Theorem \cite{Brown Halmos},
which states that the product of two scalar Toeplitz operators $T_{\Phi
}T_{\Psi}=T_{\Phi\Psi}$ ($\dim E=1$) if and only if either $\Phi$ is
co-analytic or $\Psi$ is analytic. There are many generalizations of
Brown-Halmos Theorem on different function spaces, see
\cite{Ahern,GuArchiv,Murphy}. A somewhat complicated result is obtained in the
case $\dim E>1$ in \cite{BS2}. Recently, a simple characterization of the
product of two block Toeplitz operators $T_{\Phi}T_{\Psi}$ being a block
Toeplitz operator is given in \cite{GuH}. A natural question from (\ref{h2})
is when the product of two block Hankel operators $H_{\Phi}T_{\Psi}$ is a
block Hankel operator. Indeed partial results are obtained for this question
and we will give a simple answer in Theorem \ref{hankel}.

In this paper, we mainly study the compactness of the product $H_{\Phi}%
T_{\Psi}.$ However, some compactness results hold for a (finite) sum of
(finite) products of several block Hankel operators and block Toeplitz
operators. Next we introduce the notion of an odd product.

\begin{definition}
We call the product of several block Toeplitz operators and an odd number of
block Hankel operators an odd product of block Hankel operators and block
Toeplitz operators.
\end{definition}

The proof of the next lemma shows that the product of several block Toeplitz
operators and an even number of block Hankel operators is a finite sum of
finite products of block Toeplitz operators.

\begin{lemma}
\label{product}If $A,B\in B(H_{E}^{2})$ are sums of odd products of block
Hankel operators and block Toeplitz operators, then $AB$ is a sum of products
of block Toeplitz operators. In particular, $B^{\ast}B$ is a sum of products
of block Toeplitz operators.
\end{lemma}

\begin{proof}
It is clear that we only need to prove that the product of two odd products of
block Hankel operators and block Toeplitz operators is a sum of products of
block Toeplitz operators. On a reflection, we need to prove the following
claim: the product%
\[
H_{\Phi_{1}}T_{\Psi_{1}}\cdots T_{\Psi_{m}}H_{\Phi_{2}},
\]
where $\Phi_{i},\Psi_{j}\in L_{B(E)}^{\infty},$ is a sum of finite products of
Toeplitz operators. For $m=0,$ by (\ref{t1}),%
\[
H_{\Phi_{1}}H_{\Phi_{2}}=T_{\widetilde{\Phi}_{1}\Phi_{2}}-T_{\widetilde{\Phi
}_{1}}T_{\Phi_{2}},
\]
so the claim is true. Now assume the claim holds for $m.$ Then, by
(\ref{h1}),
\begin{align*}
&  H_{\Phi_{1}}T_{\Psi_{1}}\cdots T_{\Psi_{m}}T_{\Psi_{m+1}}H_{\Phi_{2}}\\
&  =H_{\Phi_{1}}T_{\Psi_{1}}\cdots T_{\Psi_{m}}(H_{\widetilde{\Psi}_{m+1}%
\Phi_{2}}-H_{\widetilde{\Psi}_{m+1}}T_{\Phi_{2}})\\
&  =H_{\Phi_{1}}T_{\Psi_{1}}\cdots T_{\Psi_{m}}H_{\widetilde{\Psi}_{m+1}%
\Phi_{2}}-\left(  H_{\Phi_{1}}T_{\Psi_{1}}\cdots T_{\Psi_{m}}H_{\widetilde
{\Psi}_{m+1}}\right)  T_{\Phi_{2}}.
\end{align*}
By the induction hypothesis, $H_{\Phi_{1}}T_{\Psi_{1}}\cdots T_{\Psi_{m}%
}H_{\widetilde{\Psi}_{m+1}\Phi_{2}}$ and $H_{\Phi_{1}}T_{\Psi_{1}}\cdots
T_{\Psi_{m}}H_{\widetilde{\Psi}_{m+1}}$ are sums of products of block Toeplitz
operators, so the claim holds for $m+1.$
\end{proof}

Essentially, a sum of odd products of block Hankel operators and block
Toeplitz operators, denoted by $B,$ is related to a block Hankel operator in
the sense that $S_{E}^{\ast}B-BS_{E}$ is a finite rank operator. A sum of
products of block Toeplitz operators, denoted by $A,$ is related to a block
Toeplitz operator in the sense that $A-S_{E}^{\ast}AS_{E}$ is a finite rank
operator. For the study of compactness, we actually will need more general
formulas which are extensions of these observations using an automorphism of
$\mathbb{D}$.

For $z\in\mathbb{D},$ let $\varphi_{z}$ be the automorphism of $\mathbb{D}$,
\[
\varphi_{z}(w)=\frac{w-z}{1-\overline{z}w},w\in\mathbb{D}.
\]
Note that $\varphi_{0}(w)=w$ and $\widetilde{\varphi}_{z}=\overline
{\varphi_{\overline{z}}}.$ Let $k_{z}$ be the normalized reproducing kernel of
$H^{2},$%
\[
k_{z}(w)=\frac{\sqrt{1-\left\vert z\right\vert ^{2}}}{1-\overline{z}w}.
\]

We now state a few facts about block Hankel operators and block Toeplitz
operators related to $\varphi_{z}(w).$ By an abuse of notation, we will often
write $T_{\varphi_{z}I_{E}}$ and $H_{\overline{\varphi_{z}}I_{E}}$ as
$T_{\varphi_{z}}$ and $H_{\overline{\varphi_{z}}}$ respectively, where $I_{E}$
is the identity operator on $E.$

\begin{lemma}
\label{mobius}Assume $E=\mathbb{C}^{n}.$ The following holds.%
\begin{align*}
H_{\overline{\varphi_{z}}I_{E}}  &  =\sum_{i=1}^{n}k_{\overline{z}}%
e_{i}\otimes k_{z}e_{i},\\
C_{z,E}  &  :=I-T_{\varphi_{z}I_{E}}T_{\varphi_{z}I_{E}}^{\ast}\\
&  =H_{\overline{\varphi_{z}}I_{E}}^{\ast}H_{\overline{\varphi_{z}}I_{E}}%
=\sum_{i=1}^{n}k_{z}e_{i}\otimes k_{z}e_{i}.
\end{align*}

\end{lemma}

\begin{proof}
The formula for $H_{\overline{\varphi_{z}}I_{E}}$ can be verified directly. By
(\ref{t1}),
\begin{align*}
I-T_{\varphi_{z}I_{E}}T_{\varphi_{z}I_{E}}^{\ast}  &  =H_{\widetilde{\varphi
}_{z}I_{E}}H_{\overline{\varphi_{z}}I_{E}}\\
&  =\left(  \sum_{i=1}^{n}k_{z}e_{i}\otimes k_{\overline{z}}e_{i}\right)
\left(  \sum_{i=1}^{n}k_{\overline{z}}e_{i}\otimes k_{z}e_{i}\right) \\
&  =\sum_{i=1}^{n}k_{z}e_{i}\otimes k_{z}e_{i},
\end{align*}
since $\left\langle k_{\overline{z}}e_{i},k_{\overline{z}}e_{j}\right\rangle
=\delta_{i,j}.$ Note that $C_{z,E}$ is in fact the projection onto the model
space associated with the inner matrix $\varphi_{z}I_{E}.$
\end{proof}

Note that $C_{0,E}=C_{E}$ which is defined in (\ref{defce}). The above lemma
also follows from a general fact that for a two sided inner function
$\Theta\in H_{B(E)}^{\infty},$
\[
H_{\Theta^{\ast}}^{\ast}H_{\Theta^{\ast}}=I-T_{\Theta}T_{\Theta}^{\ast}%
\]
and $H_{\Theta^{\ast}}^{\ast}H_{\Theta^{\ast}}$ is the projection from
$H_{E}^{2}$ onto $H_{E}^{2}\ominus\Theta H_{E}^{2}.$ Thus $H_{\overline
{\varphi_{z}}I_{E}}$ is a partial isometry and $C_{z,E}$ is a projection. We
introduce the following two operators defined on $B(H_{E}^{2}).$ For $X\in
B(H_{E}^{2}),$ set%
\begin{align*}
\Delta_{z}(X)  &  =X-T_{\varphi_{z}I_{E}}^{\ast}XT_{\varphi_{z}I_{E}},\\
\Omega_{z}(X)  &  =XT_{\varphi_{z}I_{E}}-T_{\varphi_{\overline{z}}I_{E}}%
^{\ast}X.
\end{align*}

\begin{lemma}
Let $\Phi\in L_{B(E)}^{\infty}.$ Then
\begin{align}
T_{\varphi_{z}}^{\ast}T_{\Phi}  &  =T_{\Phi}T_{\varphi_{z}}^{\ast}%
+T_{\varphi_{z}}^{\ast}T_{\Phi}C_{z,E},\label{aa}\\
T_{\Phi}T_{\varphi_{z}}  &  =T_{\varphi_{z}}T_{\Phi}+C_{z,E}T_{\Phi}%
T_{\varphi_{z}},\label{bb}\\
H_{\Phi}T_{\varphi_{z}I_{E}}^{\ast}  &  =T_{\varphi_{\overline{z}}I_{E}%
}H_{\Phi}-T_{\varphi_{\overline{z}}I_{E}}H_{\Phi}C_{z,E}+C_{\overline{z}%
,E}H_{\Phi}T_{\varphi_{z}I_{E}}^{\ast}. \label{ccc}%
\end{align}

\end{lemma}

\begin{proof}
Since $T_{\varphi_{z}}^{\ast}T_{\Phi}T_{\varphi_{z}}=T_{\Phi},$ we have%
\begin{align*}
T_{\varphi_{z}}^{\ast}T_{\Phi}T_{\varphi_{z}}T_{\varphi_{z}}^{\ast}  &
=T_{\Phi}T_{\varphi_{z}}^{\ast},\\
T_{\varphi_{z}}^{\ast}T_{\Phi}\left(  I-C_{z,E}\right)   &  =T_{\Phi
}T_{\varphi_{z}}^{\ast},
\end{align*}
and (\ref{aa}) follows. By taking adjoint of (\ref{aa}) and noticing that the
adjoint of a block Toeplitz operator is still a block Toeplitz operator, we
have (\ref{bb}).

By Corollary \ref{cort2h2},
\[
H_{\varphi_{z}\Phi}=H_{\Phi}T_{\varphi_{z}I_{E}}=T_{\widetilde{\varphi}%
_{z}I_{E}}H_{\Phi}=T_{\varphi_{\overline{z}}I_{E}}^{\ast}H_{\Phi}.
\]
Multiplying both sides on the right by $T_{\varphi_{z}I_{E}}^{\ast}$ and on
the left by $T_{\varphi_{\overline{z}}I_{E}},$ we have%
\begin{align*}
T_{\varphi_{\overline{z}}I_{E}}T_{\varphi_{\overline{z}}I_{E}}^{\ast}H_{\Phi
}T_{\varphi_{z}I_{E}}^{\ast}  &  =T_{\varphi_{\overline{z}}I_{E}}H_{\Phi
}T_{\varphi_{z}I_{E}}T_{\varphi_{z}I_{E}}^{\ast},\\
(I-C_{\overline{z},E})H_{\Phi}T_{\varphi_{z}I_{E}}^{\ast}  &  =T_{\varphi
_{\overline{z}}I_{E}}H_{\Phi}(I-C_{z,E}).
\end{align*}
Now (\ref{ccc}) follows.
\end{proof}

By simple computations, we see that if $X$ is a sum of products of block
Toeplitz operators or a sum of odd products of block Hankel operators and
block Toeplitz operators, $\Delta_{z}(A)$ or $\Omega_{z}(B)$ is a finite rank operator.

\begin{theorem}
\label{thm1}Assume $E=\mathbb{C}^{n}.$ If $A\in B(H_{E}^{2})$ is a sum of
products of block Toeplitz operators on the vector-valued Hardy space
$H_{E}^{2}$, then $\Delta_{z}(A)$ is a finite rank operator.
\end{theorem}

\begin{proof}
Let $A=T_{\Psi_{1}}\cdots T_{\Psi_{m}},$ where $\Psi_{j}\in L_{B(E)}^{\infty
}.$ By using (\ref{bb}) iteratively and noting that $C_{z,E}$ is of finite
rank, we see that
\[
AT_{\varphi_{z}}=T_{\varphi_{z}}A+F,
\]
where $F$ is a finite rank operator. Multiplying both sides by $T_{\varphi
_{z}}^{\ast},$ we have
\begin{align*}
T_{\varphi_{z}}^{\ast}AT_{\varphi_{z}}  &  =T_{\varphi_{z}}^{\ast}%
T_{\varphi_{z}}A+T_{\varphi_{z}}^{\ast}F\\
&  =A+T_{\varphi_{z}}^{\ast}F.
\end{align*}
Thus $\Delta_{z}(A)=-T_{\varphi_{z}}^{\ast}F$ is a finite rank operator.
\end{proof}

\begin{remark}
\label{rm1}Set $z=0,$ then $\varphi_{0}=w,$ so $T_{\varphi_{z}}=S_{E}$. We see
that $\Delta_{0}(A)=A-S_{E}^{\ast}AS_{E}$ is a finite rank operator. In the
scalar case, this fact is noted in \cite{GuJFA}. The converse is also
essentially true. For example, assume $\dim E=1$ and $A\in B(H^{2})$ satisfies%
\[
A-S^{\ast}AS=\sum_{i=1}^{n}\varphi_{i}\otimes\psi_{i},
\]
where $\varphi_{i},\psi_{i}\in H^{2}.$ Then
\begin{align*}
&  \Delta_{0}(A-\sum_{i=1}^{n}H_{J\varphi_{i}}H_{J\psi_{i}}^{\ast})=\Delta
_{0}(A)-\Delta_{0}(\sum_{i=1}^{n}H_{J\varphi_{i}}H_{J\psi_{i}}^{\ast})\\
&  =\sum_{i=1}^{n}\varphi_{i}\otimes\psi_{i}-\sum_{i=1}^{n}\varphi_{i}%
\otimes\psi_{i}=0.
\end{align*}
Thus $A-\sum_{i=1}^{n}H_{J\varphi_{i}}H_{J\psi_{i}}^{\ast}=T_{\varphi},$ where
$\varphi\in L^{2}.$ This converse part is noted and used in an essential way
in \cite{GuJFA2}. The above argument is not rigorous since $H_{J\varphi_{i}%
}H_{J\psi_{i}}^{\ast}$ and $T_{\varphi}$ are potentially unbounded densely
defined operators; see \cite{GuJFA2,GuH} for a more rigorous treatment of
products of unbounded scalar or block Hankel operators and block Toeplitz operators.
\end{remark}

\begin{theorem}
\label{thm2}Assume $E=\mathbb{C}^{n}.$ If $B\in B(H_{E}^{2})$ is a sum of odd
products of block Hankel operators and block Toeplitz operators on the
vector-valued Hardy space $H_{E}^{2}$, then $\Omega_{z}(B)$ is a finite rank operator.
\end{theorem}

\begin{proof}
We first prove the case $B=T_{\Psi_{1}}H_{\Psi_{2}}T_{\Psi_{3}}.$ By
(\ref{bb}),
\[
T_{\Psi_{3}}T_{\varphi_{z}}=T_{\varphi_{z}}T_{\Psi_{3}}+F_{3},
\]
where $F_{3}$ is a finite rank operator. By Corollary \ref{cort2h2},
\[
H_{\Psi_{2}}T_{\varphi_{z}I_{E}}=T_{\widetilde{\varphi}_{z}I_{E}}H_{\Psi_{2}%
}=T_{\varphi_{\overline{z}}I_{E}}^{\ast}H_{\Psi_{2}}.
\]
By (\ref{aa}), we have%
\[
T_{\Psi_{1}}T_{\varphi_{\overline{z}}}^{\ast}=T_{\varphi_{\overline{z}}}%
^{\ast}T_{\Psi_{1}}+F_{1},
\]
where $F_{1}$ is a finite rank operator. Combining the above computations, we
see that
\[
\Omega_{z}(T_{\Psi_{1}}H_{\Psi_{2}}T_{\Psi_{3}})=T_{\Psi_{1}}H_{\Psi_{2}%
}T_{\Psi_{3}}T_{\varphi_{z}I_{E}}-T_{\varphi_{\overline{z}}I_{E}}^{\ast
}T_{\Psi_{1}}H_{\Psi_{2}}T_{\Psi_{3}}%
\]
is of finite rank.

In the general case, let $B=C_{1}\cdots C_{m},$ where each $C_{i}$ is either a
block Toeplitz operator or a block Hankel operator and an odd number of
$C_{i}^{\prime}s$ are block Hankel operators. The strategy is to start with
$BT_{\varphi_{z}}=C_{1}\cdots C_{m}T_{\varphi_{z}}$, and then to move
$T_{\varphi_{z}}$ or $T_{\varphi_{\overline{z}}}^{\ast}$ to the left across
$C_{i}^{\prime}s$ until we reach $T_{\varphi_{\overline{z}}}^{\ast}B.$ If
$C_{i}$ is a block Hankel operator, we apply either formula $C_{i}%
T_{\varphi_{z}}=T_{\varphi_{\overline{z}}}^{\ast}C_{i}$ or (\ref{ccc}) with
$\overline{z}$ in place of $z$ to get
\[
C_{i}T_{\varphi_{\overline{z}}}^{\ast}=T_{\varphi_{z}}C_{i}+F_{i},
\]
where $F_{i}$ is of finite rank. If $C_{i}$ is a block Toeplitz operator, we
apply either the formula (\ref{bb}) or (\ref{aa}) to get
\begin{align*}
C_{i}T_{\varphi_{z}}  &  =T_{\varphi_{z}}C_{i}+F_{i},\\
\text{or }C_{i}T_{\varphi_{\overline{z}}}^{\ast}  &  =T_{\varphi_{\overline
{z}}}^{\ast}C_{i}+F_{i},
\end{align*}
where $F_{i}$ is of finite rank. We conclude that $\Omega_{z}(B)$ is a finite
rank operator.
\end{proof}

The following more explicit results are crucial for us, which have their own interest.

\begin{lemma}
\label{key}Assume $E=\mathbb{C}^{n}$ and $\Phi,\Psi\in L_{B(E)}^{\infty}.$
Then
\begin{align*}
\Omega_{z}(H_{\Phi}T_{\Psi})  &  =H_{\Phi}H_{\overline{\varphi_{\overline{z}}%
}}H_{\Psi}=\sum_{i=1}^{n}H_{\Phi}\left(  k_{z}e_{i}\right)  \otimes H_{\Psi
}^{\ast}\left(  k_{\overline{z}}e_{i}\right)  ,\\
\Omega_{z}(T_{\Psi}H_{\Phi})  &  =-\sum_{i=1}^{n}H_{\widetilde{\Psi}}\left(
k_{z}e_{i}\right)  \otimes H_{\Phi}^{\ast}\left(  k_{\overline{z}}%
e_{i}\right)  .
\end{align*}

\end{lemma}

\begin{proof}
Instead using (\ref{bb}), we use (\ref{t1}) to get%
\[
T_{\Psi}T_{\varphi_{z}}=T_{\varphi_{z}\Psi}=T_{\varphi_{z}}T_{\Psi
}+H_{\overline{\varphi_{\overline{z}}}}H_{\Psi}.
\]
Then
\begin{align*}
H_{\Phi}T_{\Psi}T_{\varphi_{z}}  &  =H_{\Phi}(T_{\varphi_{z}}T_{\Psi
}+H_{\overline{\varphi_{\overline{z}}}}H_{\Psi})\\
&  =H_{\Phi}T_{\varphi_{z}}T_{\Psi}+H_{\Phi}H_{\overline{\varphi_{\overline
{z}}}}H_{\Psi}\\
&  =T_{\varphi_{\overline{z}}}^{\ast}H_{\Phi}T_{\Psi}+H_{\Phi}H_{\overline
{\varphi_{\overline{z}}}}H_{\Psi}.
\end{align*}
Thus, by Lemma \ref{mobius},
\[
\Omega_{z}(H_{\Phi}T_{\Psi})=H_{\Phi}H_{\overline{\varphi_{\overline{z}}}%
I_{E}}H_{\Psi}=\sum_{i=1}^{n}H_{\Phi}\left(  k_{z}e_{i}\right)  \otimes
H_{\Psi}^{\ast}\left(  k_{\overline{z}}e_{i}\right)  .
\]
The formula for $\Omega_{z}(T_{\Psi}H_{\Phi})$ can be proved similarly.
\end{proof}

\begin{remark}
\label{rm2}Set $z=0$ in Theorem \ref{thm2}, and we see that $\Omega
_{0}(B)=BS_{E}-S_{E}^{\ast}B$ is a finite rank operator. The converse is
essentially also true. For example, assume $\dim E=1$ and $B\in B(H^{2})$
satisfies%
\[
BS-S^{\ast}B=\sum_{i=1}^{n}\varphi_{i}\otimes\psi_{i},
\]
where $\varphi_{i},\psi_{i}\in H^{2}.$ Then, by the above lemma,
\begin{align*}
&  \Omega_{0}(B-\sum_{i=1}^{n}H_{J\varphi_{i}}T_{w\psi_{i}}^{\ast} )
=\Omega_{0}(B)-\Omega_{0}(\sum_{i=1}^{n}H_{J\varphi_{i}}T_{w\psi_{i}}^{\ast
})\\
&  =\sum_{i=1}^{n}\varphi_{i}\otimes\psi_{i}-\sum_{i=1}^{n}H_{J\varphi_{i}%
}e_{i}\otimes H_{\overline{w\psi_{i}}}^{\ast}e_{i}\\
&  =\sum_{i=1}^{n}\varphi_{i}\otimes\psi_{i}-\sum_{i=1}^{n}\varphi_{i}%
\otimes\psi_{i}=0.
\end{align*}
Thus $B-\sum_{i=1}^{n}H_{J\varphi_{i}}T_{w\psi_{i}}^{\ast}=H_{\varphi},$ where
$\varphi\in L^{2}.$ The above argument is not rigorous since $H_{J\varphi_{i}%
}T_{w\psi_{i}}^{\ast}$ and $H_{\varphi}$ are potentially unbounded densely
defined operators, see \cite{GuJFA2,GuH} for a more rigorous treatment of
products of unbounded scalar or block Hankel and Toeplitz operators.
\end{remark}

We note that a sum of odd products of block Hankel operators and block
Toeplitz operators can be written as a sum of products of block Toeplitz
operators only when it is a finite rank operator.

\begin{proposition}
Assume $E=\mathbb{C}^{n}$. If $B\in B(H_{E}^{2})$ is a sum of odd products of block Hankel operators and
block Toeplitz operators, and $B$ is also a product of Toeplitz operators,
Then $B$ is a finite rank operator.
\end{proposition}

\begin{proof}
By Remark \ref{rm1} and Remark \ref{rm2}, there exist two finite rank
operators $F_{1}$ and $F_{2}$ such that
\[
B{S_{E}}-S^{\ast}_{E}B=F_{1},B-S^{\ast}_{E}B{S_{E}}=F_{2}.
\]
Then%
\begin{align*}
B  &  =S^{\ast}_{E}B{S_{E}}+F_{2}=(BS_{E}-F_{1}){S_{E}}+F_{2},\\
B(I-{S_{E}^{2}})  &  =-F_{1}S_{E}+F_{2}.
\end{align*}
Since the range of $(I-S_{E}^{2})$ is dense in $H_{E}^{2},$ it follows that
$B$ is a finite rank operator.
\end{proof}

Let $M_{n\times r}$ denote the space of $n\times r$ complex matrices, $H$ be a
complex separable Hilbert space. For $d>0,$ set%
\[
\left(  M_{n\times r}\right)  _{d}=\left\{  C=\left[  c_{i,j}\right]  \in
M_{n\times r}:\left\vert c_{i,j}\right\vert \leq d\text{ for all }i,j\right\}
.
\]
Set $M_{n}:=M_{n\times n}$. Let $A=\left[  a_{i,j}\right]  \in M_{n}$, define
\[
\Vert A\Vert_{\infty}=\sup_{1\leq i,j\leq n}\left\vert a_{i,j}\right\vert .
\]
We need to use the fact that $\left(  M_{n}\right)  _{d}$ is a compact subset
of $M_{n}$ later on.

\begin{lemma}
\label{xy0}Let $x_{i},y_{i},z_{i}\in H.$ Let $A\in M_{r\times n}.$ If
\begin{align*}
\left[
\begin{array}
[c]{ccc}%
x_{1} & \cdots & x_{n}%
\end{array}
\right]   &  =\left[
\begin{array}
[c]{ccc}%
z_{1} & \cdots & z_{r}%
\end{array}
\right]  A,\\
\left[
\begin{array}
[c]{ccc}%
w_{1} & \cdots & w_{r}%
\end{array}
\right]   &  =\left[
\begin{array}
[c]{ccc}%
y_{1} & \cdots & y_{n}%
\end{array}
\right]  A^{\ast},
\end{align*}
then%
\[
\sum_{i=1}^{n}x_{i}\otimes y_{i}=\sum_{i=1}^{r}z_{i}\otimes w_{i}.
\]

\end{lemma}

\begin{proof}
The proof is by a direct verification. For clarity, we prove the lemma for
case $n=r=2.$ Let
\begin{align*}
\left[
\begin{array}
[c]{cc}%
x_{1} & x_{2}%
\end{array}
\right]   &  =\left[
\begin{array}
[c]{cc}%
z_{1} & z_{2}%
\end{array}
\right]  \left[
\begin{array}
[c]{cc}%
a & b\\
c & d
\end{array}
\right]  ,\\
\left[
\begin{array}
[c]{cc}%
w_{1} & w_{2}%
\end{array}
\right]   &  =\left[
\begin{array}
[c]{cc}%
y_{1} & y_{2}%
\end{array}
\right]  \left[
\begin{array}
[c]{cc}%
\overline{a} & \overline{c}\\
\overline{b} & \overline{d}%
\end{array}
\right]  ,
\end{align*}
Then%
\begin{align*}
&  x_{1}\otimes y_{1}+x_{2}\otimes y_{2}\\
&  =(az_{1}+cz_{2})\otimes y_{1}+(bz_{1}+dz_{2})\otimes y_{2}\\
&  =z_{1}\otimes\overline{a}y_{1}+z_{2}\otimes\overline{c}y_{1}+z_{1}%
\otimes\overline{b}y_{2}+z_{2}\otimes\overline{d}y_{2}\\
&  =z_{1}\otimes\left(  \overline{a}y_{1}+\overline{b}y_{2}\right)
+z_{2}\otimes\left(  \overline{c}y_{1}+\overline{d}y_{2}\right) \\
&  =z_{1}\otimes w_{1}+z_{2}\otimes w_{2}.
\end{align*}
The proof is complete.
\end{proof}

The following lemma is essentially Proposition 4 in \cite{GuZheng} with an
appropriate reformulation. Its proof helps to introduce some notations and
ideas which will be useful later on, see Lemma \ref{xy2} and the proof of
Theorem \ref{crucial}, so we include a slightly improved proof.

\begin{lemma}
\label{xy}Let $x_{i},y_{i}\in H.$ Then
\[
\sum_{i=1}^{n}x_{i}\otimes y_{i}=0
\]
if and only if there exists $A\in\left(  M_{n}\right)  _{1}$ such that
\begin{equation}
\left[
\begin{array}
[c]{ccc}%
x_{1} & \cdots & x_{n}%
\end{array}
\right]  (I-A)=0\text{ and }\left[
\begin{array}
[c]{ccc}%
y_{1} & \cdots & y_{n}%
\end{array}
\right]  A^{\ast}=0. \label{permute0}%
\end{equation}

\end{lemma}

\begin{proof}
Let $x:=\left[  x_{1},\cdots,x_{n}\right]  $, $y:=\left[  y_{1},\cdots
,y_{n}\right]  $. Then for a permutation $\sigma$ of $\{1,2,\ldots,n\},$
\[
x_{\sigma}:=\left[  x_{\sigma(1)},\cdots,x_{\sigma(n)}\right]  ,y_{\sigma
}:=\left[  y_{\sigma(1)},\cdots,y_{\sigma(n)}\right]  .
\]
We first prove that $\sum_{i=1}^{n}x_{i}\otimes y_{i}=0$ if and only if there
are a matrix $A_{0}\in\left(  M_{n}\right)  _{1}$ and a permutation $\sigma$
such that
\begin{equation}
x_{\sigma}\left(  I-A_{0}\right)  =0\quad\text{ and }\quad y_{\sigma}%
A_{0}^{\ast}=0. \label{permute}%
\end{equation}
Let $R$ be the permutation matrix such that $x=x_{\sigma}R$, $y=y_{\sigma}R$
and $A=R^{\ast}A_{0}R$. Then (\ref{permute0}) holds if and only if
(\ref{permute}) holds, since%
\begin{align*}
x(I-A)  &  =x_{\sigma}R(I-R^{\ast}A_{0}R)=x_{\sigma}\left(  I-A_{0}\right)
R,\\
yA  &  =y_{\sigma}R(R^{\ast}A_{0}^{\ast}R)=y_{\sigma}A_{0}^{\ast}R.
\end{align*}
Therefore, we need to prove (\ref{permute}).

We now prove the "if" direction of (\ref{permute}) using Lemma \ref{xy0}.
Assume (\ref{permute}) holds. Set
\begin{align*}
f  &  :=\left[  f_{1},\cdots,f_{n}\right]  =x_{\sigma}(I-A_{0}),z:=\left[
z_{1},\cdots,z_{n}\right]  =x_{\sigma}A_{0},\\
g  &  :=\left[  g_{1},\cdots,g_{n}\right]  =y_{\sigma}A_{0}^{\ast}.
\end{align*}
Since $x_{\sigma}=x_{\sigma}(I-A_{0})+x_{\sigma}A_{0}=f+z,$ we have
\begin{align}
\sum_{i=1}^{n}x_{i}\otimes y_{i}  &  =\sum_{i=1}^{n}x_{\sigma(i)}\otimes
y_{\sigma(i)}\nonumber\\
&  =\sum_{i=1}^{n}f_{i}\otimes y_{\sigma(i)}+\sum_{i=1}^{n}z_{i}\otimes
y_{\sigma(i)}\nonumber\\
&  =\sum_{i=1}^{n}f_{i}\otimes y_{\sigma(i)}+\sum_{i=1}^{n}x_{\sigma_{i}%
}\otimes g_{i}\label{xya}\\
&  =0,\nonumber
\end{align}
where the third equality follows from Lemma \ref{xy0} and the last equality
follows from (\ref{permute}).

We use induction to prove\ the "only if" direction of (\ref{permute}). Assume
$\sum_{i=1}^{n}x_{i}\otimes y_{i}=0.$ It is clear that for $n=1$, the result
is true with $A_{0}=1$ or $A_{0}=0$. Now assume the result is true for $n-1$.
Without loss of generality, assume that
\[
\max_{1\leq i\leq n}\left\Vert y_{i}\right\Vert =\left\Vert y_{j}\right\Vert
>0
\]
for some $j$. Then
\[
\sum_{i=1}^{n}\left(  x_{i}\otimes y_{i}\right)  y_{j} =\sum_{i=1}%
^{n}\left\langle y_{j},y_{i}\right\rangle x_{i} =0.
\]
Thus
\begin{equation}
x_{j}+\sum_{i\neq j}a_{i}x_{i}=0, \label{xj}%
\end{equation}
where
\[
a_{i}=\frac{\left\langle y_{j},y_{i}\right\rangle }{\left\langle y_{j}%
,y_{j}\right\rangle },\quad\left\vert a_{i}\right\vert \leq1\text{ for }i\neq
j.
\]
Now we rewrite $\sum_{i=1}^{n}x_{i}\otimes y_{i}$ as%
\[
\sum_{i=1}^{n}x_{i}\otimes y_{i}=\left(  x_{j}+\sum_{i\neq j}a_{i}%
x_{i}\right)  \otimes y_{j}+\sum_{i\neq j}x_{i}\otimes\left(  y_{i}%
-\overline{a_{i}}y_{j}\right)  .
\]

It follows from (\ref{xj}) that%
\[
\sum_{i\neq j}x_{i}\otimes\left(  y_{i}-\overline{a_{i}}y_{j}\right)  =0.
\]
By induction hypothesis, there exist $A_{1}\in\left(  M_{n-1}\right)  _{1}$
and a permutation $\omega$ of $\{1,\ldots,$ $j-1,j+1,\ldots,n\}$ such that
\begin{align}
\left[  x_{\omega(1)},\cdots,x_{\omega(j-1)},x_{\omega(j+1)},\cdots
,x_{\omega(n)}\right]  (I-A_{1})  &  =0,\label{inductioina}\\
\left(  \left[  y_{\omega(1)},\cdots,y_{\omega(j-1)},y_{\omega(j+1)}%
,\cdots,y_{\omega(n)}\right]  -z\right)  A_{1}^{\ast}  &  =0,
\label{induction2}%
\end{align}
where%
\[
z=\left[  \overline{a_{\omega(1)}}y_{j},\ldots,\overline{a_{\omega(j-1)}}%
y_{j},\overline{a_{\omega(j+1)}}y_{j},\ldots,\overline{a_{\omega(n)}}%
y_{j}\right]  .
\]
Let
\[
A_{0}=\left[
\begin{array}
[c]{cc}%
0 & 0\\
-a & A_{1}%
\end{array}
\right]  ,\quad\text{ where }\quad a=\left[  a_{\omega(1)},\ldots
,a_{\omega(j-1)},a_{\omega(j+1)},\ldots,a_{\omega(n)}\right]  ^{tr}.
\]
Take $\sigma$ to be such that $\sigma(1)=j$, $\sigma(i+1)=\omega(i)$ for
$1\leq i\leq j-1$ and $\sigma(i)=\omega(i)$ for $j+1\leq i\leq n$. It follows
from (\ref{xj}), (\ref{inductioina}) and (\ref{induction2}) that
(\ref{permute}) holds.
\end{proof}

\begin{theorem}
\label{hankel}Assume $E=\mathbb{C}^{n}$ and $\Phi,\Psi\in L_{B(E)}^{\infty}.$
Then $H_{\Phi}T_{\Psi}$ is a block Hankel operator on the vector-valued Hardy
space $H_{E}^{2}$ if and only if there exists $A\in\left(  M_{n}\right)  _{1}$
such that
\[
\Phi(z)\left[  I-A\right]  ,A\Psi(z)\in H_{B(E)}^{\infty}.
\]
In this case, $H_{\Phi}T_{\Psi}=H_{\Phi A\Psi}.$
\end{theorem}

\begin{proof}
We first prove the "if" direction. Since $\Phi-\Phi A\in H_{B(E)}^{\infty},$
we have%
\[
H_{\Phi}T_{\Psi}=H_{\Phi A}T_{\Psi}=H_{\Phi}T_{A\Psi}=H_{\Phi A\Psi},
\]
where the second equality holds since $A$ is a constant matrix and the last
equality is by (\ref{h2}).

Assume now $H_{\Phi}T_{\Psi}$ is a block Hankel operator. By Lemma \ref{key},
\[
\Omega_{0}(H_{\Phi}T_{\Psi})=\sum_{i=1}^{n}H_{\Phi}\left(  e_{i}\right)
\otimes H_{\Psi}^{\ast}\left(  e_{i}\right)  .
\]
By Lemma \ref{xy}, there exists $A\in\left(  M_{n}\right)  _{1}$ such that
\begin{align*}
\left[
\begin{array}
[c]{ccc}%
H_{\Phi}\left(  e_{1}\right)  & \cdots & H_{\Phi}\left(  e_{n}\right)
\end{array}
\right]  (I-A)  &  =0,\\
\left[
\begin{array}
[c]{ccc}%
H_{\Psi}^{\ast}\left(  e_{1}\right)  & \cdots & H_{\Psi}^{\ast}\left(
e_{n}\right)
\end{array}
\right]  A^{\ast}  &  =0.
\end{align*}
Equivalently, $\Phi(z)\left[  I-A\right]  ,A\Psi(z)\in H_{B(E)}^{\infty}.$
\end{proof}

The above theorem gives a simple and complete answer to the question when
$H_{\Phi}T_{\Psi}$ is a block Hankel operator. Some partial and more
complicated answers are obtained before \cite{GuJOT}. However, one can still
ask what $\Phi$ and $\Psi$ look like. Inspired by the answers to the question
of when $T_{\Phi}T_{\Psi}$ is a block Toeplitz operator \cite{{GuH}}, we have
the following more explicit description of $\Phi$ and $\Psi$. We recall a
lemma from \cite{GuH}.

\begin{lemma}
\label{xy1}Let $x_{i},y_{i}\in H.$ Then
\[
\sum_{i=1}^{n}x_{i}\otimes y_{i}=0
\]
if and only if there exists an invertible matrix $U\in M_{n}$ and $0\leq l\leq
n$ such that
\begin{align*}
\left[
\begin{array}
[c]{ccc}%
x_{1} & \cdots & x_{n}%
\end{array}
\right]   &  =\left[
\begin{array}
[c]{cccccc}%
u_{1} & \cdots & u_{l} & 0 & \cdots & 0
\end{array}
\right]  U\\
\left[
\begin{array}
[c]{ccc}%
y_{1} & \cdots & y_{n}%
\end{array}
\right]   &  =\left[
\begin{array}
[c]{cccccc}%
0 & \cdots & 0 & w_{l+1} & \cdots & w_{n}%
\end{array}
\right]  U^{\ast-1}%
\end{align*}
for some $u_{i}\in H$ for $1\leq i\leq l$ and $w_{j}\in H$ for $l+1\leq j\leq
n$.
\end{lemma}

Let $L_{m\times l}^{\infty}$ and $H_{m\times l}^{\infty}$ denote spaces of the
matrix-valued functions of size $m\times l$ with entries being $L^{\infty}$
functions and $H^{\infty}$ functions, respectively.

\begin{theorem}
\label{Theorem Huw} Assume $E=\mathbb{C}^{n}$ and $\Phi,\Psi\in L_{B(E)}%
^{\infty}.$ Then $H_{\Phi}T_{\Psi}$ is a block Hankel operator on the
vector-valued Hardy space $H_{E}^{2}$ if and only if there exist an invertible
matrix $D\in M_{n}$ and $0\leq l\leq n$ such that
\begin{align*}
\Phi(z)  &  =\left[
\begin{array}
[c]{cc}%
U_{1}(z) & W_{2}(z)
\end{array}
\right]  D,\\
\Psi(z)  &  =D^{-1}\left[
\begin{array}
[c]{c}%
W_{1}(z)\\
U_{2}(z)
\end{array}
\right]  ,
\end{align*}
where $U_{1}\in L_{n\times l}^{\infty},$ $U_{2}\in L_{(n-l)\times n}^{\infty
},$ $W_{1}\in H_{l\times n}^{\infty},$ $W_{2}\in H_{n\times(n-l)}^{\infty}.$
In this case $H_{\Phi}T_{\Psi}=H_{U_{1}(z)W_{1}(z)}.$
\end{theorem}

\begin{proof}
The proof is similar to the proof of Theorem \ref{hankel} by using Lemma
\ref{xy1}.
\end{proof}

\section{Proof of Theorem \ref{thm harmonic extension}}

\label{Section 3} The main purpose of this section is to prove Theorem
\ref{thm harmonic extension}. We establish some general theorems on
compactness of the product of block Toeplitz operators and block Hankel
operators, such as Theorem \ref{main}, which are useful in the next section.
These theorems will be useful in the further study of Toeplitz algebras on the
vector-valued Hardy space. We recall that for $e\in E,$ $k_{z}e\rightarrow0$
weakly as $\left\vert z\right\vert \rightarrow1.$ The following lemma provides
a property of compact operators in $B(H_{E}^{2})$, which will be used later.

\begin{lemma}
\label{ness}Assume $E=\mathbb{C}^{n}$. If $X\in B(H_{E}^{2})$ is a compact operator, then%
\[
\lim_{\left\vert z\right\vert \rightarrow1}\left\Vert \Delta_{z}(X)\right\Vert
=0\text{ and }\lim_{\left\vert z\right\vert \rightarrow1}\left\Vert \Omega
_{z}(X)\right\Vert =0\text{.}%
\]

\end{lemma}

\begin{proof}
The $\lim_{\left\vert z\right\vert \rightarrow1}\left\Vert \Delta
_{z}(X)\right\Vert =0$ in the scalar case is due to Lemma 2 in \cite{Zheng},
and the extension to the vector-valued case is immediate, as we show how to
reduce the proof of $\lim_{\left\vert z\right\vert \rightarrow1}\left\Vert
\Omega_{z}(X)\right\Vert =0$ in the vector-valued case to the scalar-valued
case below.

Assume $E=\mathbb{C}^{n}.$ If $X\in B(H_{E}^{2})$ is a compact operator, then
\[
X=\left[  X_{i,j}\right]  _{i,j=1}^{n}\text{,}%
\]
where each $X_{i,j}$ is a compact operator on $H^{2}$. Thus $\lim_{\left\vert
z\right\vert \rightarrow1}\left\Vert \Omega_{z}(X)\right\Vert =0$ if and only
if $\lim_{\left\vert z\right\vert \rightarrow1}$$\left\Vert \Omega_{z}%
(X_{i,j})\right\Vert $ =$\lim_{\left\vert z\right\vert \rightarrow1}\left\Vert
X_{i,j}T_{\varphi_{z}}-T_{\varphi_{\overline{z}}}^{\ast}X_{i,j}\right\Vert =0$
for all $i,j$.

Since the set of finite rank operators is dense in the set of compact
operators, given $\varepsilon>0,$ there exist $f_{1},\cdots,f_{n}$ and
$g_{1},\cdots,g_{n}$ in $H^{2}$ so that
\[
\left\Vert X_{i,j}-\sum_{i=1}^{n}f_{i}\otimes g_{i}\right\Vert <\varepsilon.
\]
Thus, this lemma follows once we prove it for any operator of rank one.

Note that for almost all $w\in\mathbb{T}$, as $|z|\rightarrow1,$
\begin{align*}
z+\overline{{\varphi}_{\overline{z}}}(w)  &  =z+\frac{\bar{w}-z}{1-\bar{z}%
\bar{w}}=\frac{(1-|z|^{2})\bar{w}}{1-\bar{z}\bar{w}}\rightarrow0,\\
\overline{z}+\overline{{\varphi}_{z}}(w)  &  =\overline{z}+\frac{\bar
{w}-\overline{z}}{1-z\bar{w}}=\frac{(1-|z|^{2})\bar{w}}{1-z\bar{w}}%
\rightarrow0.
\end{align*}
Let $f,g\in L^{2}.$ By Lebesgue Dominated Convergence Theorem, as
$|z|\rightarrow1$,%
\[
\left\Vert zf+\overline{{\varphi}_{\overline{z}}}f\right\Vert \rightarrow
0\text{ and }\left\Vert \overline{z}g+\overline{\varphi_{z}}g\right\Vert
\rightarrow0.
\]
If $f,g\in H^{2}$, we apply the projection $P$ to obtain, as $|z|\rightarrow
1$,%
\begin{align}
\left\Vert P\left(  zf+\overline{{\varphi}_{\overline{z}}}f\right)
\right\Vert  &  =\left\Vert zf+T_{\overline{{\varphi}_{\overline{z}}}%
}f\right\Vert \rightarrow0,\label{f1}\\
\left\Vert P\left(  \overline{z}g+\overline{\varphi_{z}}f\right)  \right\Vert
&  =\left\Vert \overline{z}g+T_{\overline{\varphi_{z}}}g\right\Vert
\rightarrow0. \label{g1}%
\end{align}
Note that
\begin{align*}
&  \left\Vert zf\otimes T_{\overline{\varphi_{z}}}g-T_{\overline{{\varphi
}_{\overline{z}}}}f\otimes\overline{z}g\right\Vert \\
&  =\left\Vert \left(  zf+T_{\overline{{\varphi}_{\overline{z}}}}f\right)
\otimes T_{\overline{\varphi_{z}}}g-T_{\overline{{\varphi}_{\overline{z}}}%
}f\otimes\left(  T_{\overline{\varphi_{z}}}g+\overline{z}g\right)
)\right\Vert .
\end{align*}
Therefore,%
\begin{align*}
\left\Vert \Omega_{z}(f\otimes g)\right\Vert  &  =\left\Vert (f\otimes
g)T_{\varphi_{z}}-T_{\varphi_{\overline{z}}}^{\ast}(f\otimes g)\right\Vert \\
&  =\left\Vert zf\otimes T_{\overline{\varphi_{z}}}g-T_{\overline{{\varphi
}_{\overline{z}}}}f\otimes\overline{z}g\right\Vert +(1-\left\vert z\right\vert
)\left\Vert \Omega_{z}(f\otimes g)\right\Vert \\
&  \leq\left\Vert zf+T_{\overline{{\varphi}_{\overline{z}}}}f\right\Vert
\left\Vert g\right\Vert +\left\Vert \overline{z}g+T_{\overline{\varphi_{z}}%
}g\right\Vert \left\Vert f\right\Vert +2(1-\left\vert z\right\vert )\left\Vert
g\right\Vert \left\Vert f\right\Vert .
\end{align*}
It follows from (\ref{f1}) and (\ref{g1}) that $\lim_{\left\vert z\right\vert
\rightarrow1}\left\Vert \Omega_{z}(f\otimes g)\right\Vert =0.$
\end{proof}

\begin{theorem}
\cite{GuJFA,GuoZ2} Let $A$ be a sum of products of scalar Toeplitz operators
on $H^{2}$ with symbol zero. Then $A$ is compact if and only if
\[
\lim_{\left\vert z\right\vert \rightarrow1}\left\Vert \Delta_{z}(A)\right\Vert
=0.
\]

\end{theorem}

The extension of the above theorem to $H_{E}^{2}$ where $\dim E<\infty$ is immediate.

\begin{corollary}
\label{nesu}Assume $E=\mathbb{C}^{n}.$ Let $A$ be a sum of product of block
Toeplitz operators on $H_{E}^{2}$ with symbol zero. Then $A$ is compact if and
only if
\[
\lim_{\left\vert z\right\vert \rightarrow1}\left\Vert \Delta_{z}(A)\right\Vert
=0.
\]

\end{corollary}

\begin{proof}
It is clear that if $A$ is a sum of product of block Toeplitz operators on
$H_{E}^{2}$ with symbol zero, then%
\[
A=\left[  A_{i,j}\right]  _{i,j=1}^{n},
\]
where each $A_{i,j}$ is a sum of products of scalar Toeplitz operators on
$H^{2}$ with symbol zero. Furthermore, $A$ is compact if and only if each
$A_{i,j}$ is compact for all $i,j$ and
\[
\lim_{\left\vert z\right\vert \rightarrow1}\left\Vert \Delta_{z}(A)\right\Vert
=0
\]
if and only if $\lim_{\left\vert z\right\vert \rightarrow1}\left\Vert
\Delta_{z}(A_{i,j})\right\Vert =0$ for all $i,j.$ Thus the result follows from
the previous theorem.
\end{proof}

It is an interesting open question whether the above result extends to the
case $\dim E=\infty.$ We note that the "only if" direction is still valid when
$\dim E=\infty.$ Next we explore a relation between $\Omega_{z}(X)$ and
$\Delta_{z}(X).$

\begin{lemma}
\label{relation}Assume $E=\mathbb{C}^{n}$. Let $X,Y\in B(H_{E}^{2}).$ Then%
\[
\Delta_{z}(XY)=XC_{\overline{z},E}Y-XT_{\varphi_{\overline{z}}}\Omega
_{z}(Y)+\Omega_{\overline{z}}(X)\Omega_{z}(Y)+\Omega_{\overline{z}%
}(X)T_{\varphi_{\overline{z}}}^{\ast}Y.
\]
In particular,%
\[
\Delta_{z}(XX^{\ast})=XC_{\overline{z},E}X^{\ast}+\left[  \Omega_{\overline
{z}}(X)T_{\varphi_{\overline{z}}}^{\ast}X^{\ast}\right]  ^{\ast}%
-\Omega_{\overline{z}}(X)\Omega_{\overline{z}}(X)^{\ast}+\Omega_{\overline{z}%
}(X)T_{\varphi_{\overline{z}}}^{\ast}X^{\ast}.
\]

\end{lemma}

\begin{proof}
By definition,%
\begin{align*}
YT_{\varphi_{z}I_{E}}  &  =\Omega_{z}(Y)+T_{\varphi_{\overline{z}}I_{E}}%
^{\ast}Y,\\
T_{\varphi_{z}I_{E}}^{\ast}X  &  =XT_{\varphi_{\overline{z}}I_{E}}%
-\Omega_{\overline{z}}(X).
\end{align*}
Hence
\begin{align*}
\Delta_{z}(XY)  &  =XY-T_{\varphi_{z}}^{\ast}XYT_{\varphi_{z}}\\
&  =XY-\left[  XT_{\varphi_{\overline{z}}}-\Omega_{\overline{z}}(X)\right]
\left[  \Omega_{z}(Y)+T_{\varphi_{\overline{z}}}^{\ast}Y\right] \\
&  =XY-XT_{\varphi_{\overline{z}}}T_{\varphi_{\overline{z}}}^{\ast
}Y-XT_{\varphi_{\overline{z}}}\Omega_{z}(Y)+\Omega_{\overline{z}}(X)\Omega
_{z}(Y)+\Omega_{\overline{z}}(X)T_{\varphi_{\overline{z}}}^{\ast}Y\\
&  =X(I-T_{\varphi_{\overline{z}}}T_{\varphi_{\overline{z}}}^{\ast
})Y-XT_{\varphi_{\overline{z}}}\Omega_{z}(Y)+\Omega_{\overline{z}}%
(X)\Omega_{z}(Y)+\Omega_{\overline{z}}(X)T_{\varphi_{\overline{z}}}^{\ast}Y.
\end{align*}
The result follows by noting $C_{\overline{z},E}=I-T_{\varphi_{\overline{z}%
}I_{E}}T_{\varphi_{\overline{z}}I_{E}}^{\ast}.$

Note that when $Y=X^{\ast},$%
\[
\Omega_{z}(X^{\ast})=-\Omega_{\overline{z}}(X)^{\ast}.
\]
The proof is complete.
\end{proof}

The following result provides a necessary and sufficient condition for a sum
of odd products of block Hankel operators and block Toeplitz operators to be compact.

\begin{theorem}
\label{main}Assume $E=\mathbb{C}^{n}.$ Let $B\in B(H_{E}^{2})$ be a sum of odd
products of block Toeplitz and Hankel operators. Then $B$ is compact if and
only if
\begin{equation}
\lim_{\left\vert z\right\vert \rightarrow1}\left\Vert \Omega_{z}(B)\right\Vert
=0\text{ and }\lim_{\left\vert z\right\vert \rightarrow1}\left\Vert
BH_{\overline{\varphi_{\overline{z}}}I_{E}}\right\Vert =0. \label{limit}%
\end{equation}

\end{theorem}

\begin{proof}
We first note that since $H_{\overline{\varphi_{z}}I_{E}}$ is a partial
isometry,%
\begin{align*}
\left\Vert BH_{\overline{\varphi_{\overline{z}}}I_{E}}\right\Vert  &
=\left\Vert BH_{\overline{\varphi_{z}}I_{E}}^{\ast}\right\Vert =\left\Vert
BH_{\overline{\varphi_{z}}I_{E}}^{\ast}H_{\overline{\varphi_{z}}I_{E}%
}H_{\overline{\varphi_{z}}I_{E}}^{\ast}\right\Vert \\
&  \leq\left\Vert BH_{\overline{\varphi_{z}}I_{E}}^{\ast}H_{\overline
{\varphi_{z}}I_{E}}\right\Vert =\left\Vert BC_{z,E}\right\Vert .
\end{align*}
On the other hand,
\[
\left\Vert BC_{z,E}\right\Vert =\left\Vert BH_{\overline{\varphi_{z}}I_{E}%
}^{\ast}H_{\overline{\varphi_{z}}I_{E}}\right\Vert \leq\left\Vert
BH_{\overline{\varphi_{z}}I_{E}}^{\ast}\right\Vert =\left\Vert BH_{\overline
{\varphi_{\overline{z}}}I_{E}}\right\Vert .
\]
Thus%
\[
\lim_{\left\vert z\right\vert \rightarrow1}\left\Vert BH_{\overline
{\varphi_{\overline{z}}}I_{E}}\right\Vert =0
\]
if and only if
\[
\lim_{\left\vert z\right\vert \rightarrow1}\left\Vert BC_{z,E}\right\Vert =0,
\]
if and only if%
\[
\lim_{\left\vert z\right\vert \rightarrow1}\left\Vert BC_{\overline{z}%
,E}\right\Vert =0,
\]
since $\left\vert z\right\vert \rightarrow1$ if and only if $\left\vert
\overline{z}\right\vert \rightarrow1.$

Assume $B$ is compact. By Lemma \ref{ness}, $\lim_{\left\vert z\right\vert
\rightarrow1}\left\Vert \Omega_{z}(B)\right\Vert =0$. Since $k_{z}%
e_{i}\rightarrow0$ weakly as $\left\vert z\right\vert \rightarrow1$ for $1\leq
i\leq n,$
\[
\lim_{\left\vert z\right\vert \rightarrow1}\left\Vert B\left(  k_{z}%
e_{i}\right)  \right\Vert =0.
\]
By Lemma \ref{mobius},%
\[
\lim_{\left\vert z\right\vert \rightarrow1}\left\Vert BH_{\overline
{\varphi_{\overline{z}}}I_{E}}\right\Vert =\lim_{\left\vert z\right\vert
\rightarrow1}\left\Vert \sum_{i=1}^{n}B\left(  k_{z}e_{i}\right)
\otimes\left(  k_{z}e_{i}\right)  \right\Vert =0.
\]

Now assume (\ref{limit}) holds. Then $\lim_{\left\vert z\right\vert
\rightarrow1}\left\Vert BC_{\overline{z},E}\right\Vert =0.$ By Lemma
\ref{relation},
\begin{align*}
\left\Vert \Delta_{z}(BB^{\ast})\right\Vert  &  \leq\left\Vert BC_{\overline
{z},E}B^{\ast}\right\Vert +2\left\Vert \Omega_{\overline{z}}(B)T_{\varphi
_{\overline{z}}}^{\ast}B^{\ast}\right\Vert +\left\Vert \Omega_{\overline{z}%
}(B)\right\Vert ^{2}\\
&  \leq\left\Vert BC_{\overline{z},E}\right\Vert ^{2}+2\left\Vert B\right\Vert
\left\Vert \Omega_{\overline{z}}(B)\right\Vert +\left\Vert \Omega
_{\overline{z}}(B)\right\Vert ^{2}.
\end{align*}
Thus $\lim_{\left\vert z\right\vert \rightarrow1}\left\Vert \Delta
_{z}(BB^{\ast})\right\Vert =0.$ By Lemma \ref{product}, $BB^{\ast}$ is a sum
of products of block Toeplitz operators. It then follows from Corollary
\ref{nesu} that $BB^{\ast}$ is compact. Therefore, $B$ is compact.
\end{proof}

It follows from the above proof that $\lim_{\left\vert z\right\vert
\rightarrow1}\left\Vert \Delta_{z}(BB^{\ast})\right\Vert =0$ is also a
necessary and sufficient condition for $B$ to be compact. However, the finite
rank operator $\Omega_{z}(B)$ in (\ref{limit}) is much simpler than the finite
rank operator $\Delta_{z}(BB^{\ast})$ for further analysis.

\begin{corollary}
\label{corlllary3.6} Assume $E=\mathbb{C}^{n}$ and $\Phi,\Psi\in
L_{B(E)}^{\infty}.$ Then $H_{\Phi}T_{\Psi}$ is compact on the vector-valued
Hardy space $H_{E}^{2}$ if and only if
\begin{align}
\lim_{\left\vert z\right\vert \rightarrow1}\left\Vert \Omega_{z}(H_{\Phi
}T_{\Psi})\right\Vert  &  =\lim_{\left\vert z\right\vert \rightarrow
1}\left\Vert \sum_{i=1}^{n}H_{\Phi}\left(  k_{z}e_{i}\right)  \otimes H_{\Psi
}^{\ast}\left(  k_{\overline{z}}e_{i}\right)  \right\Vert =0\text{ and
}\label{cond1}\\
\lim_{\left\vert z\right\vert \rightarrow1}\left\Vert H_{\Phi}T_{\Psi
}H_{\overline{\varphi_{\overline{z}}}I_{E}}\right\Vert  &  =\lim_{\left\vert
z\right\vert \rightarrow1}\left\Vert \sum_{i=1}^{n}H_{\Phi}T_{\Psi}\left(
k_{z}e_{i}\right)  \otimes\left(  k_{\overline{z}}e_{i}\right)  \right\Vert
=0. \label{cond2}%
\end{align}

\end{corollary}

\begin{proof}
The result follows from Theorem \ref{main} and Lemma \ref{key}.
\end{proof}

Next we express (\ref{cond1}) and (\ref{cond2}) in terms of the harmonic
extension of symbols $\Phi$ and $\Psi.$ Let $Trace$ denote the trace on the
trace class operators on a complex Hilbert space $H$ and $trace$ denote the
trace on $M_{n}.$

\begin{lemma}
\label{trace}Let $x_{i},y_{i}\in H.$ Let $T=\sum_{i=1}^{n}x_{i}\otimes y_{i}$
be a finite rank operator on $H$. Write%
\[
W_{x}=\left[  \left\langle x_{i},x_{j}\right\rangle \right]  _{i,j=1}%
^{n},W_{y}=\left[  \left\langle y_{i},y_{j}\right\rangle \right]  _{i,j=1}%
^{n}.
\]
Then
\[
Trace(T^{\ast}T)=trace\left(  W_{x}W_{y}\right)  .
\]

\end{lemma}

\begin{proof}
Note that%
\[
T^{\ast}T=\left(  \sum_{i=1}^{n}y_{i}\otimes x_{i}\right)  \left(  \sum
_{j=1}^{n}x_{j}\otimes y_{j}\right)  =\sum_{i=1}^{n}\sum_{j=1}^{n}\left\langle
x_{j},x_{i}\right\rangle y_{i}\otimes y_{j}.
\]
Thus%
\begin{align*}
Trace(T^{\ast}T)  &  =\sum_{i=1}^{n}\sum_{j=1}^{n}\left\langle x_{j}%
,x_{i}\right\rangle Trace\left(  y_{i}\otimes y_{j}\right) \\
&  =\sum_{i=1}^{n}\sum_{j=1}^{n}\left\langle x_{j},x_{i}\right\rangle
\left\langle y_{i},y_{j}\right\rangle \\
&  =\sum_{j=1}^{n}\left[  W_{x}W_{y}\right]  _{jj}=trace\left(  W_{x}%
W_{y}\right)  .
\end{align*}
Hence, we complete the proof of the Lemma.
\end{proof}

Now, we are ready to prove Theorem \ref{thm harmonic extension}, which is
Theorem \ref{Theorem 3.8}, as follows.

\begin{theorem}
\label{Theorem 3.8} Assume $E=\mathbb{C}^{n}$ and $\Phi,\Psi\in L_{B(E)}%
^{\infty}.$ Then $H_{\Phi}T_{\Psi}$ is compact on the vector-valued Hardy
space $H_{E}^{2}$ if and only if%
\begin{align*}
\lim_{\left\vert z\right\vert \rightarrow1}trace\left[  \left\vert \Phi
_{-}-\Phi_{-}(z)\right\vert ^{2}(z)\left\vert \widetilde{\Psi_{-}}^{\ast
}-\widetilde{\Psi_{-}}^{\ast}(\overline{z})\right\vert ^{2}(\overline
{z})\right]   &  =0\text{ and }\\
\lim_{\left\vert z\right\vert \rightarrow1}trace\left[  \left\vert \left(
\Phi_{-}\Psi_{+}\right)  _{-}-\left(  \Phi_{-}\Psi_{+}\right)  _{-}(z)+\left(
\left[  \Phi_{-}-\Phi_{-}(z)\right]  \Psi_{-}(z)\right)  \right\vert
^{2}(z)\right]   &  =0.
\end{align*}

\end{theorem}

\begin{proof}
By Lemma \ref{trace},
\begin{align*}
&  Trace\left[  \Omega_{z}(H_{\Phi}T_{\Psi})^{\ast}\Omega_{z}(H_{\Phi}T_{\Psi
})\right] \\
&  =trace\left[  W_{\Phi}W_{\Psi}\right]  =trace\left[  W_{\Phi}^{tr}W_{\Psi
}^{tr}\right]  ,
\end{align*}
where%
\begin{align*}
W_{\Phi}  &  =\left[  \left\langle H_{\Phi}\left(  k_{z}e_{i}\right)
,H_{\Phi}\left(  k_{z}e_{j}\right)  \right\rangle \right]  _{i,j=1}^{n},\\
W_{\Psi}  &  =\left[  \left\langle H_{\Psi}^{\ast}\left(  k_{\overline{z}%
}e_{i}\right)  ,H_{\Psi}^{\ast}\left(  k_{\overline{z}}e_{j}\right)
\right\rangle \right]  _{i,j=1}^{n}.
\end{align*}
Write $\Phi=\left[  \varphi_{ij}\right]  _{n\times n}$ as $\Phi=\Phi_{+}%
+\Phi_{-}$ with $\Phi_{+}=\left[  \left(  \varphi_{ij}\right)  _{+}\right]  $
and $\Phi_{-}=\left[  \left(  \varphi_{ij}\right)  _{-}\right]  .$ For a
matrix $A,$ we define $\left\vert A\right\vert ^{2}=A^{\ast}A.$ Since
$H_{\Phi}$ only depends on $\Phi_{-},$ we assume $\Phi=\Phi_{-}$ and
$\Psi=\Psi_{-}$ for the moment. Then%
\begin{align*}
\left\langle H_{\Phi}\left(  k_{z}e_{i}\right)  ,H_{\Phi}\left(  k_{z}%
e_{j}\right)  \right\rangle  &  =\sum_{m=1}^{n}\left\langle (I-P)\left(
\varphi_{mi}k_{z}\right)  ,(I-P)\left(  \varphi_{mj}k_{z}\right)
\right\rangle \\
&  =\sum_{m=1}^{n}\left\langle \left[  \varphi_{mi}-\varphi_{mi}(z)\right]
k_{z},\left[  \varphi_{mj}-\varphi_{mj}(z)\right]  k_{z}\right\rangle
\end{align*}
and
\[
W_{\Phi}^{tr}=\left\vert \Phi-\Phi(z)\right\vert ^{2}(z),W_{\Psi}%
^{tr}=\left\vert \widetilde{\Psi}^{\ast}-\widetilde{\Psi}^{\ast}(\overline
{z})\right\vert ^{2}(\overline{z}).
\]
Since for a finite rank operator $F$ with $rank(F)\leq n,$
\[
\left\Vert F\right\Vert ^{2}\leq Trace(F^{\ast}F)\leq C(n)\left\Vert
F\right\Vert ^{2},
\]
where $C(n)$ is a constant that only depends on $n.$ A similar result holds if
$F\in M_{n}$ is a matrix. Therefore, (\ref{cond1}) is equivalent to
\[
\lim_{\left\vert z\right\vert \rightarrow1}trace\left[  \left\vert \Phi
_{-}-\Phi_{-}(z)\right\vert ^{2}(z)\left\vert \widetilde{\Psi_{-}}^{\ast
}-\widetilde{\Psi_{-}}^{\ast}(\overline{z})\right\vert ^{2}(\overline
{z})\right]  =0.
\]

Similarly, by Lemma \ref{trace},%
\begin{align*}
&  Trace\left[  \left(  H_{\Phi}T_{\Psi}H_{\overline{\varphi_{\overline{z}}%
}I_{E}}\right)  ^{\ast}\left(  H_{\Phi}T_{\Psi}H_{\overline{\varphi
_{\overline{z}}}I_{E}}\right)  \right] \\
&  =trace\left[  W_{1}W_{2}\right]  ,
\end{align*}
where%
\begin{align*}
W_{1}  &  =\left[  \left\langle H_{\Phi}T_{\Psi}\left(  k_{z}e_{i}\right)
,H_{\Phi}T_{\Psi}\left(  k_{z}e_{j}\right)  \right\rangle \right]
_{i,j=1}^{n} ,\\
W_{2}  &  =\left[  \left\langle k_{\overline{z}}e_{i},k_{\overline{z}}%
e_{j}\right\rangle \right]  _{i,j=1}^{n} =I.
\end{align*}
By an abuse of notation of using unbounded block Toeplitz and Hankel
operators, we have%
\begin{align*}
H_{\Phi}T_{\Psi}\left(  k_{z}e_{i}\right)   &  =H_{\Phi}\left(  T_{\Psi_{+}%
}+T_{\Psi_{-}}\right)  \left(  k_{z}e_{i}\right) \\
&  =H_{\Phi_{-}\Psi_{+}}\left(  k_{z}e_{i}\right)  +H_{\Phi}\Psi_{-}(z)\left(
k_{z}e_{i}\right)  .
\end{align*}
Therefore,%
\[
\left(  W_{1}\right)  ^{tr}=\left\vert \left(  \Phi_{-}\Psi_{+}\right)
_{-}-\left(  \Phi_{-}\Psi_{+}\right)  _{-}(z)+\left(  \left[  \Phi_{-}%
-\Phi_{-}(z)\right]  \Psi_{-}(z)\right)  \right\vert ^{2}(z),
\]
and (\ref{cond2}) is equivalent to $\lim_{\left\vert z\right\vert
\rightarrow1} trace \left[  W_{1}\right]  =0.$
\end{proof}

We extract some facts from the above proof which are useful later on.

\begin{lemma}
\label{lemma 3.9}Assume $E=\mathbb{C}^{n}$ and $\Phi\in L_{B(E)}^{\infty}.$
Then%
\begin{align*}
&  Trace(\left(  H_{\Phi}H_{\overline{\varphi_{\overline{z}}}I_{E}}\right)
^{\ast}\left(  H_{\Phi}H_{\overline{\varphi_{\overline{z}}}I_{E}}\right)  )\\
&  =Trace(\left(  H_{\overline{\varphi_{\overline{z}}}I_{E}}H_{\Phi}\right)
^{\ast}\left(  H_{\overline{\varphi_{\overline{z}}}I_{E}}H_{\Phi}\right)  )\\
&  =trace\left\vert \Phi_{-}-\Phi_{-}(z)\right\vert ^{2}(z).
\end{align*}
Consequently, $\lim_{\left\vert z\right\vert \rightarrow1}\left\Vert
H_{\overline{\varphi_{\overline{z}}}I_{E}}H_{\Phi}\right\Vert =0$ if and only
if $\lim_{\left\vert z\right\vert \rightarrow1}\left\Vert H_{\Phi}%
H_{\overline{\varphi_{\overline{z}}}I_{E}}\right\Vert =0.$
\end{lemma}

\begin{proof}
Note that%
\begin{align*}
H_{\Phi}H_{\overline{\varphi_{\overline{z}}}I_{E}}  &  =\sum_{i=1}^{n}H_{\Phi
}\left(  k_{z}e_{i}\right)  \otimes\left(  k_{\overline{z}}e_{i}\right)  ,\\
H_{\overline{\varphi_{\overline{z}}}I_{E}}H_{\Phi}  &  =\sum_{i=1}^{n}\left(
k_{z}e_{i}\right)  \otimes H_{\Phi}^{\ast}\left(  k_{\overline{z}}%
e_{i}\right)  .
\end{align*}
By the proof of Theorem \ref{Theorem 3.8},%
\[
Trace(\left(  H_{\Phi}H_{\overline{\varphi_{\overline{z}}}I_{E}}\right)
^{\ast}\left(  H_{\Phi}H_{\overline{\varphi_{\overline{z}}}I_{E}}\right)
)=trace\left\vert \Phi-\Phi(z)\right\vert ^{2}(z)
\]
and $Trace(\left(  H_{\overline{\varphi_{\overline{z}}}I_{E}}H_{\Phi}\right)
^{\ast}\left(  H_{\overline{\varphi_{\overline{z}}}I_{E}}H_{\Phi}\right)
)=traceW$, where%
\[
W=\left[  \left\langle H_{\Phi}^{\ast}\left(  k_{\overline{z}}e_{i}\right)
,H_{\Phi}^{\ast}\left(  k_{\overline{z}}e_{j}\right)  \right\rangle \right]
_{i,j=1}^{n} .
\]
Note that (for the moment assume $\Phi=\Phi_{-}$)%
\begin{align*}
\left\langle H_{\Phi}^{\ast}\left(  k_{\overline{z}}e_{i}\right)  ,H_{\Phi
}^{\ast}\left(  k_{\overline{z}}e_{j}\right)  \right\rangle  &  =\sum
_{m=1}^{n}\left\langle (I-P)\left(  \overline{\varphi_{im}(\overline{w}%
)}k_{\overline{z}}\right)  ,(I-P)\left(  \overline{\varphi_{jm}(\overline{w}%
)}k_{\overline{z}}\right)  \right\rangle \\
&  =\sum_{m=1}^{n}\left\langle (I-P)\left(  \varphi_{jm}(w)k_{z}\right)
,(I-P)\left(  \varphi_{im}(w)k_{z}\right)  \right\rangle ,
\end{align*}
where we use the fact for $f,g\in L^{2},$%
\[
\left\langle \overline{f(\overline{w})},\overline{g(\overline{w}%
)}\right\rangle =\left\langle g(w),f(w)\right\rangle \text{ and }%
\overline{k_{\overline{z}}(\overline{w})}=k_{\overline{z}}(w).
\]
Therefore, $W=\left\vert \Phi-\Phi(z)\right\vert ^{2}(z).$ The result follows
from the proof of previous theorem.
\end{proof}

\section{Proof of Theorem \ref{maina}.}

\label{Section 4}

The section is devoted to the proof of Theorem \ref{maina}, which is long and
complicated. For a better presentation, we motivate the definitions of several
expressions related to the study of the compactness of $H_{\Phi}T_{\Psi}$; we
divide the proof of Theorem \ref{maina} into several lemmas and theorems; and
introduce various tools along the way. Careful attention is also devoted to
finding relatively simple notations to convey our ideas effectively.

Next we introduce a quantity which expresses (\ref{cond1}) and (\ref{cond2})
together in a way for us to use the theory of maximal ideal space and
localization techniques. We first give heuristical arguments which also
introduce notations. Heuristically, (\ref{cond1}) is equivalent to $\left\Vert
H_{\Phi}\left(  k_{z}e_{i}\right)  \right\Vert $ goes to zero for several
different $i$ and $\left\Vert H_{\Psi}^{\ast}\left(  k_{\overline{z}}%
e_{i}\right)  \right\Vert $ goes to zero for other different $i.$ Inspired by
Lemma \ref{xy} and Theorem \ref{hankel} which deals with the case when the
finite rank operator $\Omega_{z}(H_{\Phi}T_{\Psi})=0,$ heuristically, to show
$\Omega_{z}(H_{\Phi}T_{\Psi})$ goes to zero is equivalent to show%
\begin{align*}
&  \left[
\begin{array}
[c]{ccc}%
H_{\Phi}\left(  k_{z}e_{1}\right)  & \cdots & H_{\Phi}\left(  k_{z}%
e_{n}\right)
\end{array}
\right]  (I-A)\\
&  \text{ and }\left[
\begin{array}
[c]{ccc}%
H_{\Psi}^{\ast}\left(  k_{\overline{z}}e_{1}\right)  & \cdots & H_{\Psi}%
^{\ast}\left(  k_{\overline{z}}e_{n}\right)
\end{array}
\right]  \text{ }A^{\ast}%
\end{align*}
go to zero for some $A\in\left(  M_{n}\right)  _{1}.$ Note that%
\begin{align}
&  \left[
\begin{array}
[c]{ccc}%
H_{\Phi}\left(  k_{z}e_{1}\right)  & \cdots & H_{\Phi}\left(  k_{z}%
e_{n}\right)
\end{array}
\right]  (I-A)\label{onema}\\
&  =\left[
\begin{array}
[c]{ccc}%
H_{\Phi(I-A)}\left(  k_{z}e_{1}\right)  & \cdots & H_{\Phi(I-A)}\left(
k_{z}e_{n}\right)
\end{array}
\right]  ,\nonumber\\
&  \left[
\begin{array}
[c]{ccc}%
H_{\Psi}^{\ast}\left(  k_{\overline{z}}e_{1}\right)  & \cdots & H_{\Psi}%
^{\ast}\left(  k_{\overline{z}}e_{n}\right)
\end{array}
\right]  \text{ }A^{\ast}\nonumber\\
&  =\left[
\begin{array}
[c]{ccc}%
H_{A\Psi}^{\ast}\left(  k_{\overline{z}}e_{1}\right)  & \cdots & H_{A\Psi
}^{\ast}\left(  k_{\overline{z}}e_{n}\right)
\end{array}
\right]  .\nonumber
\end{align}
For notational simplicity, we also use the more compact form of $\Omega
_{z}(H_{\Phi}T_{\Psi})=H_{\Phi}H_{\overline{\varphi_{\overline{z}}}I_{E}%
}H_{\Psi}.$ These considerations lead to the definition of%
\[
\Gamma_{1}(z):=\inf\left\{  \left\Vert H_{\Phi(I-A)}H_{\overline
{\varphi_{\overline{z}}}I_{E}}\right\Vert +\left\Vert H_{\overline
{\varphi_{\overline{z}}}I_{E}}H_{A\Psi}\right\Vert :A\in\left(  M_{n}\right)
_{2^{2n}}\right\}  .
\]
The reason we need $\left(  M_{n}\right)  _{2^{2n}}$ instead of $\left(
M_{n}\right)  _{1}$ as in Lemma \ref{xy} is due to Lemma \ref{xy2} below.

We also need an equivalent definition of $\Gamma_{1}(z).$ Note that
\begin{align*}
H_{\Phi(I-A)}H_{\overline{\varphi_{\overline{z}}}I_{E}}  &  =\sum_{i=1}%
^{n}H_{\Phi(I-A)}\left(  k_{z}e_{i}\right)  \otimes\left(  k_{\overline{z}%
}e_{i}\right)  ,\\
H_{\overline{\varphi_{\overline{z}}}I_{E}}H_{A\Psi}  &  =\sum_{i=1}^{n}\left(
k_{z}e_{i}\right)  \otimes H_{A\Psi}^{\ast}\left(  k_{\overline{z}}%
e_{i}\right)  .
\end{align*}
Since $\left\{  k_{\overline{z}}e_{i},1\leq i\leq n\right\}  $ is an
orthonormal set, if we define%
\[
\widehat{\Gamma}_{1}(z)=\inf\left\{  \sum_{i=1}^{n}\left\Vert H_{\Phi
(I-A)}\left(  k_{z}e_{i}\right)  \right\Vert +\sum_{i=1}^{n}\left\Vert
H_{A\Psi}^{\ast}\left(  k_{\overline{z}}e_{i}\right)  \right\Vert :A\in\left(
M_{n}\right)  _{2^{2n}}\right\}  ,
\]
then there exist positive constants $C_{1}$ and $C_{2}$ that only depend on
$n$ such that
\begin{equation}
C_{1}\Gamma_{1}(z)\leq\widehat{\Gamma}_{1}(z)\leq C_{2}\Gamma_{1}(z).
\label{gammahat}%
\end{equation}
Note that for $A\in M_{n},$
\begin{align*}
&  \left\Vert H_{\Phi}H_{\overline{\varphi_{\overline{z}}}I_{E}}H_{\Psi
}\right\Vert \\
&  =\left\Vert H_{\Phi(I-A)+\Phi A}H_{\overline{\varphi_{\overline{z}}}I_{E}%
}H_{\Psi}\right\Vert \\
&  =\left\Vert H_{\Phi(I-A)}H_{\overline{\varphi_{\overline{z}}}I_{E}}H_{\Psi
}+H_{\Phi}H_{\overline{\varphi_{\overline{z}}}I_{E}}H_{A\Psi}\right\Vert \\
&  \leq\left\Vert H_{\Phi(I-A)}H_{\overline{\varphi_{\overline{z}}}I_{E}%
}\right\Vert \left\Vert H_{\Psi}\right\Vert +\left\Vert H_{\Phi}\right\Vert
\left\Vert H_{\overline{\varphi_{\overline{z}}}I_{E}}H_{A\Psi}\right\Vert .
\end{align*}
Since $A\in M_{n}$ is arbitrary, we have
\begin{equation}
\left\Vert H_{\Phi}H_{\overline{\varphi_{\overline{z}}}I_{E}}H_{\Psi
}\right\Vert \leq C\Gamma_{1}(z), \label{gamma1}%
\end{equation}
where $C=\max\{\left\Vert \Phi\right\Vert _{\infty},\left\Vert \Psi\right\Vert
_{\infty}\}.$

The above inequality immediately implies the following lemma.

\begin{lemma}
\label{easy}Notation as above. If $\lim_{\left\vert z\right\vert \rightarrow
1}\Gamma_{1}(z)=0$, then $\lim_{\left\vert z\right\vert \rightarrow
1}\left\Vert \Omega_{z}(H_{\Phi}T_{\Psi})\right\Vert =0$.
\end{lemma}

Next we are going to show if $\lim_{\left\vert z\right\vert \rightarrow
1}\left\Vert \Omega_{z}(H_{\Phi}T_{\Psi})\right\Vert =0,$ then $\lim
_{\left\vert z\right\vert \rightarrow1}\Gamma_{1}(z)=0.$ The proof is rather
difficult. We begin by proving several lemmas. Let
\[
\left(  L_{r\times r}\right)  _{1}=\left\{  C\in\left(  M_{r}\right)
_{1}:C\text{ is lower triangular with diagonals all equal to }1\right\}  .
\]

\begin{lemma}
Let $x_{i}\in H$, where $H$ is a vector space. If $\dim Span\{x_{1}%
,\cdots,x_{n}\}=r$ where $0<r<n,$ then there exist $L\in\left(  L_{(n-r)\times
(n-r)}\right)  _{1},$ $B\in\left(  M_{r\times(n-r)}\right)  _{1}$ and a
permutation $\sigma$ of $\{1,2,\ldots,n\}$ such that%
\begin{align}
Span\{x_{\sigma(n-r+1)},\cdots,x_{\sigma(n)}\}  &  =Span\{x_{1},\cdots
,x_{n}\},\label{dim}\\
\left[
\begin{array}
[c]{ccc}%
x_{\sigma(1)} & \cdots & x_{\sigma(n)}%
\end{array}
\right]  \left[
\begin{array}
[c]{c}%
L\\
B
\end{array}
\right]   &  =0 . \label{bmatrix}%
\end{align}

\end{lemma}

\begin{proof}
We prove the lemma by using induction on $n.$ When $n=2,$ it is easy to see
the lemma holds. Since $r<n,$ there exist $a_{i}\in\mathbb{C}$, not all zero,
such that $\sum_{i=1}^{n}a_{i}x_{i}=0.$ Let
\[
\max_{1\leq i\leq n}\left\vert a_{i}\right\vert =\left\vert a_{j}\right\vert
>0.
\]
Then%
\begin{equation}
x_{j}+\sum_{i\neq j}b_{i}x_{i}=0, \label{xja}%
\end{equation}
where $\left\vert b_{i}\right\vert =\left\vert a_{i}/a_{j}\right\vert \leq1.$
If $r=n-1,$ then (\ref{bmatrix}) holds with $\sigma(1)=j,$ $L=1,$ and $B$
being the column vectors consisting of $b_{i}^{\prime}s.$ If $r<n-1,$ we apply
the induction hypothesis to
\[
Span\{x_{1},\cdots,x_{j-1},x_{j+1},\cdots,x_{n}\}
\]
to obtain $L_{1}\in\left(  L_{(n-1-r)\times(n-1-r)}\right)  _{1}$, matrix
$B_{1}\in\left(  M_{r\times(n-r-1)}\right)  _{1}$ and a permutation $\omega$
of $\{1,\ldots,$ $j-1,j+1,\ldots,n\}$ such that%
\begin{equation}
\left[
\begin{array}
[c]{cccccc}%
x_{\omega(1)} & \cdots & x_{\omega(j-1)} & x_{\omega(j+1)} & \cdots &
x_{\omega(n)}%
\end{array}
\right]  \left[
\begin{array}
[c]{c}%
L_{1}\\
B_{1}%
\end{array}
\right]  =0. \label{b1}%
\end{equation}
Let $\sigma$ be such that $\sigma(1)=j,$ $\sigma(i+1)=\omega(i)$ for $1\leq
i\leq j-1$ and $\sigma(i)=\omega(i)$ for $j+1\leq i\leq n$. Set%
\[
b:=\left[  b_{\sigma(2)},\cdots,b_{\sigma(n)}\right]  ^{tr}.
\]
Then, it is easy to see that
\[
\left[
\begin{array}
[c]{cc}%
1 & 0\\
b &
\begin{array}
[c]{c}%
L_{1}\\
B_{1}%
\end{array}
\end{array}
\right]  =\left[
\begin{array}
[c]{c}%
L\\
B
\end{array}
\right]
\]
for some $L\in\left(  L_{(n-r)\times(n-r)}\right)  _{1},$ $B\in\left(
M_{r\times(n-r)}\right)  _{1}.$ Now (\ref{dim}) and (\ref{bmatrix}) follow
from (\ref{xja}) and (\ref{b1}).
\end{proof}

\begin{lemma}
\label{xy2}Let $x_{i}\in H$, where $H$ is a vector space. If $\dim
Span\{x_{1},\cdots,x_{n}\}=r$ where $0<r<n,$ then there exist $B\in\left(
M_{r\times(n-r)}\right)  _{2^{2n}}$ and a permutation $\sigma$ of
$\{1,2,\ldots,n\}$ such that%
\[
\left[
\begin{array}
[c]{ccc}%
x_{\sigma(1)} & \cdots & x_{\sigma(n-r)}%
\end{array}
\right]  =\left[
\begin{array}
[c]{ccc}%
x_{\sigma(n-r+1)} & \cdots & x_{\sigma(n)}%
\end{array}
\right]  B.
\]

\end{lemma}

\begin{proof}
We can write (\ref{bmatrix}) as%
\[
\left[
\begin{array}
[c]{ccc}%
x_{\sigma(1)} & \cdots & x_{\sigma(n-r)}%
\end{array}
\right]  =-\left[
\begin{array}
[c]{ccc}%
x_{\sigma(n-r+1)} & \cdots & x_{\sigma(n)}%
\end{array}
\right]  BL^{-1}.
\]
By using back substitution, we can see that each entry of $L^{-1}$ is bounded
by $2^{n}.$ Hence each entry of $BL^{-1}$ is bounded by $n2^{n}$ which is less
than $2^{2n}.$
\end{proof}

Let $\varphi\in L^{\infty}.$ Without loss of generality we may assume that
$\left\|  \varphi\right\|  _{\infty}<1$. As in \cite{Volberg}, there is a
unimodular function $u$ in $\varphi+H^{\infty}$ such that $T_{u}$ is invertible.

\begin{lemma}
\cite[Lemma 3]{Zheng}\label{Lemma3 Zheng} If $T_{u}$ is invertible, then there
is a constant $C_{u}>0$ such that
\[
\left\|  H_{u} k_{z}\right\|  ^{2} \leqslant\left(  1-|u(z)|^{2}\right)
\leqslant C_{u}\left\|  H_{u} k_{z}\right\|  ^{2}
\]
for all $z$ in $\mathbb{D}$.
\end{lemma}

In order to state the local conditions, it is necessary to introduce some
notation for the maximal ideal space. The space of nonzero multiplicative
linear functionals of the Douglas algebra $B_{1}$ will be denoted by
$M(B_{1})$. A function $f \in H^{\infty} $ may be thought of as a continuous
function on $M\left(  H^{\infty}+C\right)  $. For a function $F $ on the unit
disk $\mathbb{D}$ and $m \in M(H^{\infty}+ C) $, we say
\[
\lim_{z \to m} F(z) = 0
\]
if for every net $\{z_{\alpha}\} \subset\mathbb{D}$ converging to $m $,
\[
\lim_{z_{\alpha}\to m} F(z_{\alpha}) = 0.
\]
For $m$ in $M\left(  H^{\infty}\right)  $, we can identify $m$ with a
multiplicative functional on $H^{\infty}$. Furthermore, the Gleason-Whitney
theorem \cite{Douglas} says that $m$ uniquely extends to a bounded positive
linear functional $l_{m}$ on $L^{\infty}$. By the Riesz representation
theorem, there exists a measure $\mu_{m}$ called the representing measure for
$m$ with support $s(m)$, which is a subset of $M\left(  L^{\infty}\right)  $
such that
\begin{equation}
\label{support}l_{m}(f) = \int_{s(m)} f \, d\mu_{m}.
\end{equation}
For more details, see page 181 in \cite{Hoffman}. A subset of $M\left(
L^{\infty}\right)  $ is called a support set if it is the (closed) support of
the representing measure for a functional in $M\left(  H^{\infty}+C\right)  $.
The closed support of this measure $\mu_{m}$ is called the support set of $m$
and denoted by $s(m)$.

The next lemma, stated below, comes from Lemmas 2.5 and 2.6 in \cite{Gorkin},
with the proof utilizing Lemma \ref{Lemma3 Zheng}, and interprets the local
condition in an elementary way. And the extension to a vector-valued or a
matrix-valued function is straightforward by applying the scalar version
entrywise, so we will refer to this lemma in the vector-valued or
matrix-valued case.

\begin{lemma}
\label{1999Gorkin lemma} Let $f$ be in $L^{\infty}$, $m\in M\left(  H^{\infty
}+C\right)  $, and let $s(m)$ be the support set for $m$. Then the following
are equivalent: \newline(1) $f|_{s(m)}\in H^{\infty}|_{s(m)}$. \newline(2)
$\underline{\lim}_{z\rightarrow m}\left\Vert H_{f}k_{z}\right\Vert =0$.
\newline(3) ${\lim}_{z\rightarrow m}\left\Vert H_{f}k_{z}\right\Vert =0$.
\end{lemma}

\begin{theorem}
\label{crucial}Notation as above. Then $\lim_{\left\vert z\right\vert
\rightarrow1}\left\Vert \Omega_{z}(H_{\Phi}T_{\Psi})\right\Vert =0$ if and
only if $\lim_{\left\vert z\right\vert \rightarrow1}\Gamma_{1}(z)=0$.
\end{theorem}

\begin{proof}
By Lemma \ref{easy}, we need to prove the "only if" direction. Assume
$\lim_{\left\vert z\right\vert \rightarrow1}\left\Vert \Omega_{z}(H_{\Phi
}T_{\Psi})\right\Vert =0$. Suppose $\lim_{\left\vert z\right\vert
\rightarrow1}\Gamma_{1}(z)=0$ does not hold. That is, there are $\delta>0$ and
a net $\{z\}\subset\mathbb{D}$ accumulating a point in $\mathbb{T}$ such that%

\[
\widehat{\Gamma}_{1}(z)\geq\delta.
\]
We will get a contradiction. We may assume that the net $\{z\}$ converges to
some nontrivial point $m\in M\left(  H^{\infty}+C\right)  $. Then by
Carleson's Corona Theorem \cite{Garnett}, $\lim_{\left\vert z\right\vert
\rightarrow1}\left\Vert \Omega_{z}(H_{\Phi}T_{\Psi})\right\Vert =0$ implies
that $\lim_{z\rightarrow m}\left\Vert \Omega_{z}(H_{\Phi}T_{\Psi})\right\Vert
=0.$

We claim that $\lim_{z\rightarrow m}\left\Vert \Omega_{z}(H_{\Phi}T_{\Psi
})\right\Vert =0$ implies there exists $A\in\left(  M_{n}\right)  _{2^{2n}}$
such that
\begin{equation}
\lim_{z\rightarrow m}\sum_{i=1}^{n}\left\Vert H_{\Phi(I-A)}\left(  k_{z}%
e_{i}\right)  \right\Vert =\lim_{z\rightarrow m}\sum_{i=1}^{n}\left\Vert
H_{A\Psi}^{\ast}\left(  k_{\overline{z}}e_{i}\right)  \right\Vert =0.
\label{claim}%
\end{equation}
But%
\[
\sum_{i=1}^{n}\left(  \left\Vert H_{\Phi(I-A)}\left(  k_{z}e_{i}\right)
\right\Vert +\left\Vert H_{A\Psi}^{\ast}\left(  k_{\overline{z}}e_{i}\right)
\right\Vert \right)  \geq\widehat{\Gamma}_{1}(z).
\]
Hence%
\[
\lim_{z\rightarrow m}\left(  \sum_{i=1}^{n}\left(  \left\Vert H_{\Phi
(I-A)}\left(  k_{z}e_{i}\right)  \right\Vert +\left\Vert H_{A\Psi}^{\ast
}\left(  k_{\overline{z}}e_{i}\right)  \right\Vert \right)  \right)
\geq\delta.
\]
This is a contradiction.

Now we prove the claim. Let $M_{n\times1}$ denote the space of $n\times1$
complex matrices. $L_{n\times1}^{\infty}$ and $H_{n\times1}^{\infty}$ denote
spaces of size $n\times1$ with entries being $L^{\infty}$ and $H^{\infty}$
functions, respectively. Let $(L_{n\times1}^{\infty}|_{s(m)})/(H_{n\times
1}^{\infty}|_{s(m)})$ be the quotient vector space. For a function $f$ in
$L_{n\times1}^{\infty}$, let $\left\{  f\right\}  _{m}$ denote the element in
$(L_{n\times1}^{\infty}|_{s(m)})/(H_{n\times1}^{\infty}|_{s(m)})$ which
contains $f$. Write the matrix $\Phi$ by its columns as
\begin{align*}
\Phi &  =\left[
\begin{array}
[c]{ccc}%
\Phi_{1} & \cdots & \Phi_{n}%
\end{array}
\right]  ,\\
\lbrack\Phi]_{m} &  =\left[
\begin{array}
[c]{ccc}%
\left\{  \Phi_{1}\right\}  _{m} & \cdots & \left\{  \Phi_{n}\right\}  _{m}%
\end{array}
\right]  ,
\end{align*}
where $\Phi_{i}\in L_{n\times1}^{\infty}$ for all $i=1,\cdots,n$.

Suppose the dimension of the space spanned by $\left\{  \Phi_{1}\right\}
_{m},\cdots,\left\{  \Phi_{n}\right\}  _{m}$ is $N\leq n$. We assume $0<N<n.$
The proof for the case $N=0$ or $N=n$ is similar. By Lemma \ref{xy2}, up to a
permutation, we may assume that $\left\{  \left\{  \Phi_{1}\right\}
_{m},\cdots,\left\{  \Phi_{N}\right\}  _{m}\right\}  $ is a basis such that
\[
\left[
\begin{array}
[c]{ccc}%
\left\{  \Phi_{1}\right\}  _{m} & \cdots & \left\{  \Phi_{n}\right\}  _{m}%
\end{array}
\right]  =\left[
\begin{array}
[c]{ccc}%
\left\{  \Phi_{1}\right\}  _{m} & \cdots & \left\{  \Phi_{N}\right\}  _{m}%
\end{array}
\right]  B
\]
with $B=\left[  b_{ij}\right]  $ and $\left\vert b_{ij}\right\vert \leq2^{2n}%
$. The case of permutation can be treated as in the proof of Lemma \ref{xy1}.
Let $A$ be the $n\times n$ matrix $\left[
\begin{array}
[c]{c}%
B\\
0
\end{array}
\right]  $. Then $\Phi(I-A)|_{s(m)}\in$ $H_{B(E)}^{\infty}|_{s(m)}$. By Lemma
\ref{1999Gorkin lemma},%

\[
\lim_{z\rightarrow m}\sum_{i=1}^{n}\left\Vert H_{\Phi(I-A)}\left(  k_{z}%
e_{i}\right)  \right\Vert =0.
\]

On the other hand, by Lemma \ref{key}, and similar calculations as in
(\ref{xya}),
\begin{align*}
\Omega_{z}(H_{\Phi}T_{\Psi})  &  =\sum_{i=1}^{n}H_{\Phi}\left(  k_{z}%
e_{i}\right)  \otimes H_{\Psi}^{\ast}\left(  k_{\overline{z}}e_{i}\right) \\
&  =\sum_{i=1}^{n}H_{\Phi(I-A)}\left(  k_{z}e_{i}\right)  \otimes H_{\Psi
}^{\ast}\left(  k_{\overline{z}}e_{i}\right)  +\sum_{i=1}^{n}H_{\Phi}\left(
k_{z}e_{i}\right)  \otimes H_{A\Psi}^{\ast}\left(  k_{\overline{z}}%
e_{i}\right) \\
&  =\sum_{i=1}^{n}H_{\Phi(I-A)}\left(  k_{z}e_{i}\right)  \otimes H_{\Psi
}^{\ast}\left(  k_{\overline{z}}e_{i}\right)  +\sum_{j=1}^{N}H_{\Phi_{(N)}%
}\left(  k_{z}e_{j}\right)  \otimes H_{B\Psi}^{\ast}\left(  k_{\overline{z}%
}e_{j}\right)  ,
\end{align*}
where $\Phi_{(N)}=\left(  \Phi_{1},\cdots,\Phi_{N}\right)  $. Since
\[
\lim_{z\rightarrow m}\left\Vert \sum_{i=1}^{n}H_{\Phi(I-A)}\left(  k_{z}%
e_{i}\right)  \otimes H_{\Psi}^{\ast}\left(  k_{\overline{z}}e_{i}\right)
\right\Vert =0
\]
and $\lim_{z\rightarrow m}\left\Vert \Omega_{z}(H_{\Phi}T_{\Psi})\right\Vert
=0,$ we have
\begin{equation}
\lim_{z\rightarrow m}\left\Vert \sum_{j=1}^{N}H_{\Phi_{(N)}}\left(  k_{z}%
e_{j}\right)  \otimes H_{B\Psi}^{\ast}\left(  k_{\overline{z}}e_{j}\right)
\right\Vert =0. \label{phin}%
\end{equation}

We are going to show that%
\[
\lim_{z\rightarrow m}\sum_{j=1}^{N}\left\Vert H_{B\Psi}^{\ast}\left(
k_{\overline{z}}e_{j}\right)  \right\Vert =0.
\]
Suppose that this is not true. We may assume that%
\[
\lim_{z\rightarrow m}\left\Vert H_{B\Psi}^{\ast}\left(  k_{\overline{z}}%
e_{1}\right)  \right\Vert >0.
\]
Let $a_{j}(z)=\left\langle H_{B\Psi}^{\ast}\left(  k_{\overline{z}}%
e_{1}\right)  ,H_{B\Psi}^{\ast}\left(  k_{\overline{z}}e_{j}\right)
\right\rangle $. Note that $\left\vert a_{j}(z)\right\vert \leq\left\Vert
B\Psi\right\Vert _{\infty}^{2}\leq C(n)\left\Vert \Psi\right\Vert _{\infty
}^{2},$ where $C(n)$ is a constant that only depends on $n$. We may assume
that $a_{j}(z)$ converges to $\alpha_{j}\in\mathbb{C}$ as $z$ goes to $m$.
Note that $\alpha_{1}=\lim_{z\rightarrow m}\left\Vert H_{B\Psi}^{\ast}\left(
k_{\overline{z}}e_{1}\right)  \right\Vert >0$. It follows from (\ref{phin})
that
\begin{align*}
0  &  =\lim_{z\rightarrow m}\left\Vert \sum_{j=1}^{N}\left(  H_{\Phi_{(N)}%
}\left(  k_{z}e_{j}\right)  \otimes H_{B\Psi}^{\ast}\left(  k_{\overline{z}%
}e_{j}\right)  \right)  H_{B\Psi}^{\ast}\left(  k_{\overline{z}}e_{1}\right)
\right\Vert \\
&  =\lim_{z\rightarrow m}\left\Vert \sum_{j=1}^{N}a_{j}(z)H_{\Phi_{(N)}%
}\left(  k_{z}e_{j}\right)  \right\Vert =\lim_{z\rightarrow m}\left\Vert
\sum_{j=1}^{N}\alpha_{j}H_{\Phi_{(N)}}\left(  k_{z}e_{j}\right)  \right\Vert .
\end{align*}
By Lemma \ref{1999Gorkin lemma}, $\sum_{j=1}^{N}\alpha_{j}\Phi_{j}\in
H_{n\times1}^{\infty}|_{s(m)}$. This contradicts the fact that $\left\{
\left\{  \Phi_{1}\right\}  _{m},\cdots,\left\{  \Phi_{N}\right\}
_{m}\right\}  $ is a basis. Therefore,%
\[
\lim_{z\rightarrow m}\sum_{j=1}^{N}\left\Vert H_{B\Psi}^{\ast}\left(
k_{\overline{z}}e_{j}\right)  \right\Vert =\lim_{z\rightarrow m}\sum_{i=1}%
^{n}\left\Vert H_{A\Psi}^{\ast}\left(  k_{\overline{z}}e_{i}\right)
\right\Vert =0.
\]
The proof is complete.
\end{proof}

In the above theorem, we analyze $\lim_{\left\vert z\right\vert \rightarrow
1}\left\Vert \Omega_{z}(H_{\Phi}T_{\Psi})\right\Vert $ to show how the terms
in the finite rank operator $\Omega_{z}(H_{\Phi}T_{\Psi})$ interact with each
other. In order to prove Theorem \ref{maina}, we need a way to combine
(\ref{cond1}) and (\ref{cond2}). Since $\left\{  k_{\overline{z}}e_{i},1\leq
i\leq n\right\}  $ is an orthonormal set, (\ref{cond2}) holds if and only if
\begin{equation}
\lim_{\left\vert z\right\vert \rightarrow1}\left\Vert H_{\Phi}T_{\Psi}\left(
k_{z}e_{i}\right)  \right\Vert =0,1\leq i\leq n. \label{cond2a}%
\end{equation}
However, in order to combine the above condition with $\Gamma_{1}(z),$
heuristically, we need to get rid of $T_{\Psi}$ in the expression $H_{\Phi
}T_{\Psi}\left(  k_{z}e_{i}\right)  .$ A way of doing this is suggested by the
following two lemmas.

\begin{lemma}
\label{lemma 3.7} Notation as above. Then
\[
H_{\Phi}T_{\Psi}H_{\overline{\varphi_{\overline{z}}}}=H_{\Phi(I-A)}T_{\Psi
}H_{\overline{\varphi_{\overline{z}}}}+H_{\Phi A\Psi}H_{\overline
{\varphi_{\overline{z}}}}-T_{\widetilde{\Phi}}H_{A\Psi}H_{\overline
{\varphi_{\overline{z}}}}.
\]

\end{lemma}

\begin{proof}
Note that%
\begin{align*}
H_{\Phi}T_{\Psi}  &  =H_{\Phi(I-A)+\Phi A}T_{\Psi}\\
&  =H_{\Phi(I-A)}T_{\Psi}+H_{\Phi}T_{A\Psi}\\
&  =H_{\Phi(I-A)}T_{\Psi}+H_{\Phi A\Psi}-T_{\widetilde{\Phi}}H_{A\Psi}.
\end{align*}
The result follows by multiplying $H_{\overline{\varphi_{\overline{z}}}I_{E}}$
on the right.
\end{proof}

To handle the term $H_{\Phi(I-A)}T_{\Psi}H_{\overline{\varphi_{\overline{z}}}%
},$ we also need the following lemma, which in the scalar case is Lemma 17 in
\cite{GuoZ1}.

\begin{lemma}
\label{lemma 3.8} Notation as above. If $\lim_{\left\vert z\right\vert
\rightarrow1}\left\Vert H_{\Phi}H_{\overline{\varphi_{\overline{z}}}%
}\right\Vert =0,$ then $\lim_{\left\vert z\right\vert \rightarrow1}\left\Vert
H_{\Phi}T_{\Psi}H_{\overline{\varphi_{\overline{z}}}}\right\Vert =0.$
\end{lemma}

\begin{proof}
We use (\ref{cond2a}) to reduce the vector-valued case to the scalar case.
Write
\[
\Phi=\left[  \varphi_{ij}\right]  _{n\times n},\Psi=\left[  \psi_{ij}\right]
_{n\times n}.
\]
By (\ref{cond2a}), $\lim_{\left\vert z\right\vert \rightarrow1}\left\Vert
H_{\Phi}H_{\overline{\varphi_{\overline{z}}}}\right\Vert =0$ if and only if
\begin{equation}
\lim_{\left\vert z\right\vert \rightarrow1}\left\Vert H_{\varphi_{ij}}%
k_{z}\right\Vert =0\text{ for all }i,j. \label{assump}%
\end{equation}
And $\lim_{\left\vert z\right\vert \rightarrow1}\left\Vert H_{\Phi}T_{\Psi
}H_{\overline{\varphi_{\overline{z}}}}\right\Vert =0$ if and only if
\begin{equation}
\lim_{\left\vert z\right\vert \rightarrow1}\left\Vert \sum_{m=1}^{n}%
H_{\varphi_{im}}T_{\psi_{mj}}k_{z}\right\Vert =0\text{ for all }i,j.
\label{conclu}%
\end{equation}
By Lemma 17 in \cite{GuoZ1}, (\ref{assump}) implies that $\lim_{\left\vert
z\right\vert \rightarrow1}\left\Vert H_{\varphi_{im}}T_{\psi_{mj}}%
k_{z}\right\Vert =0$ for all $i,j,m.$ Therefore, (\ref{conclu}) holds.
\end{proof}

Now we combine (\ref{cond1}) and (\ref{cond2}) together to define%
\[
\Gamma_{2}(z):=\inf\left\{  \left\Vert H_{\Phi(I-A)}H_{\overline
{\varphi_{\overline{z}}}I_{E}}\right\Vert +\left\Vert H_{\overline
{\varphi_{\overline{z}}}I_{E}}H_{A\Psi}\right\Vert +\left\Vert H_{\Phi A\Psi
}H_{\overline{\varphi_{\overline{z}}}I_{E}}\right\Vert :A\in\left(
M_{n}\right)  _{2^{2n}}\right\}  .
\]
By (\ref{gammahat}), an equivalent expression for $\Gamma_{2}(z)$ is
\[
\widehat{\Gamma}_{2}(z):=\inf\left\{  \sum_{i=1}^{n}\left(  \left\Vert
H_{\Phi(I-A)}\left(  k_{z}e_{i}\right)  \right\Vert +\left\Vert H_{A\Psi
}^{\ast}\left(  k_{\overline{z}}e_{i}\right)  \right\Vert +\left\Vert H_{\Phi
A\Psi}\left(  k_{z}e_{i}\right)  \right\Vert \right)  :A\in\left(
M_{n}\right)  _{2^{2n}}\right\}  .
\]

The following result is crucial for the proof of Theorem \ref{maina}.

\begin{theorem}
\label{crucial2}Assume $E=\mathbb{C}^{n}$ and $\Phi,\Psi\in L_{B(E)}^{\infty
}.$ Then $\lim_{\left\vert z\right\vert \rightarrow1}\Gamma_{2}(z)=0$ if and
only if
\begin{equation}
\lim_{\left\vert z\right\vert \rightarrow1}\left\Vert \Omega_{z}(H_{\Phi
}T_{\Psi})\right\Vert =0\text{ and }\lim_{\left\vert z\right\vert
\rightarrow1}\left\Vert H_{\Phi}T_{\Psi}H_{\overline{\varphi_{\overline{z}}%
}I_{E}}\right\Vert =0. \label{cond2b}%
\end{equation}

\end{theorem}

\begin{proof}
Assume (\ref{cond2b}) holds. Suppose $\lim_{\left\vert z\right\vert
\rightarrow1}\Gamma_{2}(z)=0$ does not hold. That is, there are $\delta>0$ and
a net $\{z\}\subset\mathbb{D}$ accumulating a point in $\mathbb{T}$ such that%

\[
\widehat{\Gamma}_{2}(z)\geq\delta.
\]
We will get a contradiction. We may assume that the net $\{z\}$ converges to
some nontrivial point $m\in M\left(  H^{\infty}+C\right)  $.

By the proof of Theorem \ref{crucial} and (\ref{claim}), there exists
$A\in\left(  M_{n}\right)  _{2^{2n}}$ such that
\begin{align}
\lim_{z\rightarrow m}\sum_{i=1}^{n}\left\Vert H_{\Phi(I-A)}\left(  k_{z}%
e_{i}\right)  \right\Vert  &  =0,\label{claim1}\\
\lim_{z\rightarrow m}\sum_{i=1}^{n}\left\Vert H_{A\Psi}^{\ast}\left(
k_{\overline{z}}e_{i}\right)  \right\Vert  &  =0.\label{claim2}%
\end{align}

By Lemma \ref{lemma 3.7},
\begin{equation}
\left\Vert H_{\Phi A\Psi}H_{\overline{\varphi_{\overline{z}}}}\right\Vert
\leq\left\Vert H_{\Phi(I-A)}T_{\Psi}H_{\overline{\varphi_{\overline{z}}}%
}\right\Vert +\left\Vert T_{\widetilde{\Phi}}H_{A\Psi}H_{\overline
{\varphi_{\overline{z}}}}\right\Vert +\left\Vert H_{\Phi}T_{\Psi}%
H_{\overline{\varphi_{\overline{z}}}}\right\Vert . \label{37}%
\end{equation}
Note that (\ref{claim2}) is equivalent to $\lim_{z\rightarrow m}\left\Vert
H_{\overline{\varphi_{\overline{z}}}I_{E}}H_{A\Psi}\right\Vert =0$. By Lemma
\ref{lemma 3.9},
\[
\lim_{z\rightarrow m}\left\Vert H_{A\Psi}H_{\overline{\varphi_{\overline{z}}%
}I_{E}}\right\Vert =0\text{ and }\lim_{z\rightarrow m}\left\Vert
T_{\widetilde{\Phi}}H_{A\Psi}H_{\overline{\varphi_{\overline{z}}}I_{E}%
}\right\Vert =0.
\]

By Lemma \ref{lemma 3.8}, with $\Phi$ replaced by $\Phi(I-A)$, (\ref{claim1})
implies that $\lim_{z\rightarrow m}\left\Vert H_{\Phi(I-A)}T_{\Psi
}H_{\overline{\varphi_{\overline{z}}}}\right\Vert =0.$ Therefore, by
(\ref{37}),
\begin{equation}
\lim_{z\rightarrow m}\left\Vert H_{\Phi A\Psi}H_{\overline{\varphi
_{\overline{z}}}}\right\Vert =0. \label{claim3}%
\end{equation}

Note that%
\[
\sum_{i=1}^{n}\left(  \left\Vert H_{\Phi(I-A)}\left(  k_{z}e_{i}\right)
\right\Vert +\left\Vert H_{A\Psi}^{\ast}\left(  k_{\overline{z}}e_{i}\right)
\right\Vert +\left\Vert H_{\Phi A\Psi}\left(  k_{z}e_{i}\right)  \right\Vert
\right)  \geq\widehat{\Gamma}_{2}(z),
\]
so
\[
\lim_{z\rightarrow m}\left(  \sum_{i=1}^{n}\left(  \left\Vert H_{\Phi
(I-A)}\left(  k_{z}e_{i}\right)  \right\Vert +\left\Vert H_{A\Psi}^{\ast
}\left(  k_{\overline{z}}e_{i}\right)  \right\Vert +\left\Vert H_{\Phi A\Psi
}\left(  k_{z}e_{i}\right)  \right\Vert \right)  \right)  \geq\delta,
\]
which is a contradiction to the combination of (\ref{claim1}), (\ref{claim2})
and (\ref{claim3}).

Assume now $\lim_{\left\vert z\right\vert \rightarrow1}\Gamma_{2}(z)=0.$ We
will show (\ref{cond2b}) holds. The proof is similar to the above proof by
reversing the steps. For clarity, we include the details. We are going to show
that for each $m\in M\left(  H^{\infty}+C\right)  ,$ if the net $\{z\}$
converges to $m,$ then
\[
\lim_{z\rightarrow m}\left\Vert \Omega_{z}(H_{\Phi}T_{\Psi})\right\Vert
=0\text{ and }\lim_{z\rightarrow m}\left\Vert H_{\Phi}T_{\Psi}H_{\overline
{\varphi_{\overline{z}}}I_{E}}\right\Vert =0.
\]
By $\lim_{\left\vert z\right\vert \rightarrow1}\Gamma_{2}(z)=0$, for each $z$
is the net $\{z\},$ there exists $A_{z}\in\left(  M_{n}\right)  _{2^{2n}}$
such that
\[
\lim_{z\rightarrow m}\left(  \sum_{i=1}^{n}\left(  \left\Vert H_{\Phi
(I-A_{z})}\left(  k_{z}e_{i}\right)  \right\Vert +\left\Vert H_{A_{z}\Psi
}^{\ast}\left(  k_{\overline{z}}e_{i}\right)  \right\Vert +\left\Vert H_{\Phi
A_{z}\Psi}\left(  k_{z}e_{i}\right)  \right\Vert \right)  \right)  =0.
\]
Since $\left(  M_{n}\right)  _{2^{2n}}$ is a compact set, we may assume
$A_{z}$ converges to $A_{m},$ then%
\begin{equation}
\lim_{z\rightarrow m}\left(  \sum_{i=1}^{n}\left(  \left\Vert H_{\Phi
(I-A_{m})}\left(  k_{z}e_{i}\right)  \right\Vert +\left\Vert H_{A_{m}\Psi
}^{\ast}\left(  k_{\overline{z}}e_{i}\right)  \right\Vert +\left\Vert H_{\Phi
A_{m}\Psi}\left(  k_{z}e_{i}\right)  \right\Vert \right)  \right)
=0.\label{claim4}%
\end{equation}
By virtue of (\ref{gamma1}), $\lim_{z\rightarrow m}\left\Vert \Omega
_{z}(H_{\Phi}T_{\Psi})\right\Vert =0.$ By Lemma \ref{lemma 3.7},
\begin{equation}
\left\Vert H_{\Phi}T_{\Psi}H_{\overline{\varphi_{\overline{z}}}I_{E}%
}\right\Vert \leq\left\Vert H_{\Phi(I-A_{m})}T_{\Psi}H_{\overline
{\varphi_{\overline{z}}}I_{E}}\right\Vert +\left\Vert H_{\Phi A_{m}\Psi
}H_{\overline{\varphi_{\overline{z}}}I_{E}}\right\Vert +\left\Vert
T_{\widetilde{\Phi}}\right\Vert \left\Vert H_{A_{m}\Psi}H_{\overline
{\varphi_{\overline{z}}}I_{E}}\right\Vert .\label{37a}%
\end{equation}
By (\ref{claim4}),
\[
\lim_{z\rightarrow m}\sum_{i=1}^{n}\left\Vert H_{\Phi(I-A_{m})}\left(
k_{z}e_{i}\right)  \right\Vert =0
\]
implies
\[
\lim_{z\rightarrow m}\left\Vert H_{\Phi(I-A_{m})}H_{\overline{\varphi
_{\overline{z}}}I_{E}}\right\Vert =0.
\]
By Lemma \ref{lemma 3.8}, the above limit implies
\[
\lim_{z\rightarrow m}\left\Vert H_{\Phi(I-A_{m})}T_{\Psi}H_{\overline
{\varphi_{\overline{z}}}I_{E}}\right\Vert =0.
\]
By (\ref{claim4}),
\[
\lim_{z\rightarrow m}\sum_{i=1}^{n}\left\Vert H_{A_{m}\Psi}^{\ast}\left(
k_{\overline{z}}e_{i}\right)  \right\Vert =0
\]
implies
\[
\lim_{z\rightarrow m}\left\Vert H_{\overline{\varphi_{\overline{z}}}I_{E}%
}H_{A_{m}\Psi}\right\Vert =0.
\]
By Lemma \ref{lemma 3.9}, the above limit implies
\[
\lim_{z\rightarrow m}\left\Vert H_{A_{m}\Psi}H_{\overline{\varphi
_{\overline{z}}}I_{E}}\right\Vert =0.
\]
By (\ref{claim4}),
\[
\lim_{z\rightarrow m}\sum_{i=1}^{n}\left\Vert H_{\Phi A_{m}\Psi}\left(
k_{z}e_{i}\right)  \right\Vert =0
\]
implies
\[
\lim_{z\rightarrow m}\left\Vert H_{\Phi A_{m}\Psi}H_{\overline{\varphi
_{\overline{z}}}I_{E}}\right\Vert =0.
\]
Now it follows from (\ref{claim4}) and (\ref{37a}) that $\lim_{z\rightarrow
m}\left\Vert H_{\Phi}T_{\Psi}H_{\overline{\varphi_{\overline{z}}}I_{E}%
}\right\Vert =0.$
\end{proof}

For the proof of Theorem \ref{maina}, we need more lemmas. The following
lemma, first precisely described by D. Sarason, was later proven by Gorkin and
Zheng in \cite{Gorkin}.

\begin{lemma}
\cite[Lemma 1.3]{Gorkin}\label{1999Gorkin lemma1.3} Let $A_{\alpha}$ be a
family of Douglas algebras, the closed subalgebras of $L^{\infty}$ containing
$H^{\infty}$. Then, $M\left(  \bigcap A_{\alpha}\right)  =\overline{\bigcup
M\left(  A_{\alpha}\right)  }$.
\end{lemma}

The famous Chang-Marshall Theorem \cite{Chang, Marshall} asserts that any
Douglas algebra $B_{1}$ is generated by $H^{\infty}$ together with the complex
conjugates of a set of inner functions. This theorem also shows that the
maximal ideal space $M(B_{1})$ completely determines $B_{1}$.

\begin{theorem}
\label{Chang-Marshall} If $B_{1}$ is a Douglas algebra, then there is a set
$B_{2}$ of inner functions in $H^{\infty}$ such that
\[
B_{1}=[H^{\infty},\overline{B_{2}} ].
\]

\end{theorem}

The following lemma is a consequence of Theorem \ref{Chang-Marshall} and
(\ref{support}).

\begin{lemma}
\cite[Lemma 1.5]{Gorkin}\label{1999Gorkin lemma1.5} Let $m\in M\left(
H^{\infty}+C\right)  $, and let $s(m)$ be the support set for $m$. Then $m\in
M\left(  H^{\infty}[f]\right)  $ if and only if $f|_{s(m)}\in H^{\infty
}|_{s(m)}$.
\end{lemma}

We conclude this section by proving several equivalent conditions for the
compactness of $H_{\Phi}T_{\Psi}$, and get Theorem \ref{maina}.

\begin{theorem}
\label{main theorem TFAE} Assume $E=\mathbb{C}^{n}$ and $\Phi,\Psi\in
L_{B(E)}^{\infty}.$ Then the following statements are equivalent.\newline$(i)$
$H_{\Phi}T_{\Psi}$ is compact on the vector-valued Hardy space $H_{E}^{2}%
$.\newline$(ii)$ $\lim_{\left\vert z\right\vert \rightarrow1}\left\Vert
\Omega_{z}(H_{\Phi}T_{\Psi})\right\Vert =0$ and $\lim_{\left\vert z\right\vert
\rightarrow1}\left\Vert H_{\Phi}T_{\Psi}H_{\overline{\varphi_{\overline{z}}%
}I_{E}}\right\Vert =0$.\newline$(iii)$ $\lim_{\left\vert z\right\vert
\rightarrow1}\Gamma_{2}(z)=0$.\newline$(iv)$ For each $m\in M(H^{\infty}+C),$
there exists $A\in\left(  M_{n}\right)  _{2^{2n}}$ such that%
\[
\Phi(I-A)|_{s(m)},A\Psi|_{s(m)},\Phi A\Psi|_{s(m)}\in H_{n\times n}^{\infty
}|_{s(m)}.
\]
$(v)$ The following relation holds:
\begin{equation}
\bigcap_{A\in\left(  M_{n}\right)  _{2^{2n}}}H^{\infty}\left[  \Phi
(I-A),A\Psi,\Phi A\Psi\right]  \subset H_{n\times n}^{\infty}+C_{n\times
n}.\label{intersection}%
\end{equation}

\end{theorem}

\begin{proof}
(i)$\Longleftrightarrow$(ii). It is proved in Corollary \ref{corlllary3.6}.

(ii)$\Longleftrightarrow$(iii). It is proved in Theorem \ref{crucial2}.

(iii)$\Longrightarrow$(iv). We are going to show that for each $m\in M\left(
H^{\infty}+C\right)  $, there exists a matrix $A_{m}\in\left(  M_{n}\right)
_{2^{2n}}$ such that%
\[
\left\{  \Phi\left(  I-A_{m}\right)  \right\}  _{m}=0,\left\{  A_{m}%
\Psi\right\}  _{m}=0,\left\{  \Phi A_{m}\Psi\right\}  _{m}=0,
\]
where $\left\{  \Phi\left(  I-A_{m}\right)  \right\}  _{m}$ is an element in
the quotient vector space $(L_{n\times n}^{\infty}|_{s(m)})/(H_{n\times
n}^{\infty}|_{s(m)}).$ Let $\left\{  z\right\}  $ be a net in $\mathbb{D}$
converging to $m$. By (iii), there exists a matrix $A_{z}\in\left(
M_{n}\right)  _{2^{2n}}$ such that%

\[
\lim_{z\rightarrow m}\sum_{i=1}^{n}\left(  \left\Vert H_{\Phi(I-A_{z})}\left(
k_{z}e_{i}\right)  \right\Vert +\left\Vert H_{A_{z}\Psi}^{\ast}\left(
k_{\overline{z}}e_{i}\right)  \right\Vert +\left\Vert H_{\Phi A_{z}\Psi
}\left(  k_{z}e_{i}\right)  \right\Vert \right)  =0.
\]
Since $\left(  M_{n}\right)  _{2^{2n}}$ is compact, we may assume that $A_{z}$
converges to $A_{m}$. Hence%
\[
\lim_{z\rightarrow m}\left(  \sum_{i=1}^{n}\left(  \left\Vert H_{\Phi
(I-A_{m})}\left(  k_{z}e_{i}\right)  \right\Vert +\left\Vert H_{A_{m}\Psi
}^{\ast}\left(  k_{\overline{z}}e_{i}\right)  \right\Vert +\left\Vert H_{\Phi
A_{m}\Psi}\left(  k_{z}e_{i}\right)  \right\Vert \right)  \right)  =0.
\]
By Lemmas \ref{lemma 3.9} and Lemmas \ref{1999Gorkin lemma}, we have%
\[
\left\{  \Phi\left(  I-A_{m}\right)  \right\}  _{m}=0,\left\{  A_{m}%
\Psi\right\}  _{m}=0,\left\{  \Phi A_{m}\Psi\right\}  _{m}=0.
\]

(iv) $\Longrightarrow$ (iii). Suppose (iii) does not hold. There exists
$\delta>0$ and a net $z$ in $\mathbb{D}$ converging to some $m\in M\left(
H^{\infty}+C\right)  $ such that%
\[
\widehat{\Gamma}_{2}(z)\geq\delta.
\]
By Statement (iv), for this $m\in M\left(  H^{\infty}+C\right)  $, there
exists a matrix $A_{m}\in\left(  M_{n}\right)  _{2^{2n}}$ such that%
\[
\left\{  \Phi\left(  I-A_{m}\right)  \right\}  _{m}=0,\left\{  A_{m}%
\Psi\right\}  _{m}=0,\left\{  \Phi A_{m}\Psi\right\}  _{m}=0.
\]

Then, by Lemma \ref{lemma 3.9} and Lemma \ref{1999Gorkin lemma},
\[
\lim_{z\rightarrow m}\left(  \sum_{i=1}^{n}\left(  \left\Vert H_{\Phi
(I-A_{m})}\left(  k_{z}e_{i}\right)  \right\Vert +\left\Vert H_{A_{m}\Psi
}^{\ast}\left(  k_{\overline{z}}e_{i}\right)  \right\Vert +\left\Vert H_{\Phi
A_{m}\Psi}\left(  k_{z}e_{i}\right)  \right\Vert \right)  \right)  =0.
\]
This is a contradiction since%
\[
\sum_{i=1}^{n}\left(  \left\Vert H_{\Phi(I-A_{m})}\left(  k_{z}e_{i}\right)
\right\Vert +\left\Vert H_{A_{m}\Psi}^{\ast}\left(  k_{\overline{z}}%
e_{i}\right)  \right\Vert +\left\Vert H_{\Phi A_{m}\Psi}\left(  k_{z}%
e_{i}\right)  \right\Vert \right)  \geq\widehat{\Gamma}_{2}(z)>\delta.
\]

(iv) $\Longrightarrow$ (v). By Theorem \ref{Chang-Marshall}, we only need to
show that%
\[
M\left(  H_{n\times n}^{\infty}+C_{n\times n}\right)  \subset M\left(
\bigcap_{A\in\left(  M_{n}\right)  _{2^{2n}}}H^{\infty}\left[  \Phi
(I-A),A\Psi,\Phi A\Psi\right]  \right)  .
\]
By Lemma \ref{1999Gorkin lemma1.5}, Statement (iv) states exactly that%
\begin{equation}
M\left(  H_{n\times n}^{\infty}+C_{n\times n}\right)  \subset\bigcup
_{A\in\left(  M_{n}\right)  _{2^{2n}}}M\left(  H^{\infty}\left[
\Phi(I-A),A\Psi,\Phi A\Psi\right]  \right)  .\label{open}%
\end{equation}
By Lemma \ref{1999Gorkin lemma1.3},%
\begin{equation}
M\left(  \bigcap_{A\in\left(  M_{n}\right)  _{2^{2n}}}H^{\infty}\left[
\Phi(I-A),A\Psi,\Phi A\Psi\right]  \right)  =\overline{\bigcup_{A\in\left(
M_{n}\right)  _{2^{2n}}}M\left(  H^{\infty}\left[  \Phi(I-A),A\Psi,\Phi
A\Psi\right]  \right)  }.\label{close}%
\end{equation}
Hence%
\[
M\left(  H_{n\times n}^{\infty}+C_{n\times n}\right)  \subset M\left(
\bigcap_{A\in\left(  M_{n}\right)  _{2^{2n}}}H^{\infty}\left[  \Phi
(I-A),A\Psi,\Phi A\Psi\right]  \right)  .
\]

(v) $\Longrightarrow$ (iv). By (\ref{open}) and (\ref{close}), we need to
show
\[
\bigcup_{A\in\left(  M_{n}\right)  _{2^{2n}}}M\left(  H^{\infty}\left[
\Phi(I-A),A\Psi,\Phi A\Psi\right]  \right)  \text{ is closed.}%
\]
In other words, let $m_{\alpha}\rightarrow m$, where
\begin{equation}
m_{\alpha}\in M\left(  H^{\infty}\left[  \Phi(I-A_{\alpha}),A_{\alpha}%
\Psi,\Phi A_{\alpha}\Psi\right]  \right)  \text{ with }A_{\alpha}\in\left(
M_{n}\right)  _{2^{2n}}.\label{alpha}%
\end{equation}
We shall show that $m\in M\left(  H^{\infty}\left[  \Phi(I-A),A\Psi,\Phi
A\Psi\right]  \right)  $ for some $A\in\left(  M_{n}\right)  _{2^{2n}}.$ 

By Lemmas \ref{1999Gorkin lemma} and \ref{1999Gorkin lemma1.5}, it is
sufficient to show that
\begin{equation}
\lim_{z\rightarrow m}\sum_{i=1}^{n}\left\Vert H_{\Phi(I-A)}\left(  k_{z}%
e_{i}\right)  \right\Vert =0,\label{limit1}%
\end{equation}%
\begin{equation}
\lim_{z\rightarrow m}\sum_{i=1}^{n}\left\Vert H_{A\Psi}\left(  k_{z}%
e_{i}\right)  \right\Vert =0,\label{limit2}%
\end{equation}
and
\begin{equation}
\lim_{z\rightarrow m}\sum_{i=1}^{n}\left\Vert H_{\Phi A\Psi}\left(  k_{z}%
e_{i}\right)  \right\Vert =0.\label{limit4}%
\end{equation}
We only prove (\ref{limit1}). A similar argument can be used to get
(\ref{limit2}) and (\ref{limit4}).

Without loss of generality, we will assume $\left\Vert \Phi\right\Vert
_{\infty}\leq1$ and $\left\Vert \Psi\right\Vert _{\infty}\leq1.$ Since
$A_{\alpha}\in\left(  M_{n}\right)  _{2^{2n}},$ it can be assumed that
$A_{\alpha}\rightarrow A$. Clearly, $A\in\left(  M_{n}\right)  _{2^{2n}}.$ By
(\ref{alpha}),
\[
\lim_{z\rightarrow m_{\alpha}}\sum_{i=1}^{n}\left\Vert H_{\Phi(I-A_{\alpha}%
)}\left(  k_{z}e_{i}\right)  \right\Vert =0.
\]
Since $\left\Vert \Phi\right\Vert _{\infty}\leq1$ and $A\in\left(
M_{n}\right)  _{2^{2n}},$ there exists a constant $\lambda>0$ such that
\[
\Vert\lambda\Phi(I-A)\Vert_{\infty}<1.
\]
Thus%
\[
\operatorname{dist}_{L^{\infty}}\left(  \lambda\Phi(I-A),H_{n\times n}%
^{\infty}\right)  \leq\Vert\lambda\Phi(I-A)\Vert_{\infty}<1.
\]
As a consequence of Adamian-Arov-Krein Theorem \cite{Garnett}, there exist
unimodular functions $u_{ij}$ $(i=1,2,\cdots,n,j=1,2,\cdots,n)$ such that
$u=\left[  u_{ij}\right]_{n\times n}  \in\Phi(I-A)+H_{n\times n}^{\infty}$. Applying Lemma \ref{Lemma3 Zheng} entrywise yields%

\begin{equation}
\sum_{i=1}^{n}\left\Vert H_{\lambda\Phi(I-A)}\left(  k_{z}e_{i}\right)
\right\Vert \leq\sum_{i=1}^{n}\sum_{j=1}^{n}\left(  1-\left\vert
u_{ij}(z)\right\vert ^{2}\right)  ^{1/2}\leq C\sum_{i=1}^{n}\left\Vert
H_{\lambda\Phi(I-A)}\left(  k_{z}e_{i}\right)  \right\Vert ,\label{unimodular}%
\end{equation}
where each $u_{ij}(z)$ denotes the value of the harmonic extension of $u_{ij}$
at $z$.

Note that%
\begin{align*}
\sum_{i=1}^{n}\left\Vert H_{\lambda\Phi(I-A)}\left(  k_{z}e_{i}\right)
\right\Vert  & \leq\sum_{i=1}^{n}\left\Vert H_{\lambda\Phi(I-A_{\alpha}%
)}\left(  k_{z}e_{i}\right)  \right\Vert +\sum_{i=1}^{n}\left\Vert
H_{\lambda\Phi(A-A_{\alpha})}\left(  k_{z}e_{i}\right)  \right\Vert \\
& \leq\sum_{i=1}^{n}\left\Vert H_{\lambda\Phi(I-A_{\alpha})}\left(  k_{z}%
e_{i}\right)  \right\Vert +n\lambda\Vert A-A_{\alpha}\Vert_{\infty}.
\end{align*}
Hence
\begin{align*}
\limsup_{z\rightarrow m_{\alpha}}\sum_{i=1}^{n}\left\Vert H_{\lambda\Phi
(I-A)}\left(  k_{z}e_{i}\right)  \right\Vert  & \leq\limsup_{z\rightarrow
m_{\alpha}}\sum_{i=1}^{n}\left\Vert H_{\lambda\Phi(I-A_{\alpha})}\left(
k_{z}e_{i}\right)  \right\Vert +n\lambda\Vert A-A_{\alpha}\Vert_{\infty}\\
& =n\lambda\Vert A-A_{\alpha}\Vert_{\infty}.
\end{align*}
Together with (\ref{unimodular}), we get
\begin{align*}
\limsup_{z\rightarrow m_{\alpha}}\sum_{i=1}^{n}\sum_{j=1}^{n}\left(
1-\left\vert u_{ij}(z)\right\vert ^{2}\right)  ^{1/2}  & \leq C\limsup
_{z\rightarrow m_{\alpha}}\sum_{i=1}^{n}\left\Vert H_{\lambda\Phi(I-A)}\left(
k_{z}e_{i}\right)  \right\Vert \\
& \leq Cn\lambda\Vert A-A_{\alpha}\Vert_{\infty}.
\end{align*}
Since each $u_{ij}$ is continuous on $M\left(  H^{\infty}\right)  $, we have
\[
\sum_{i=1}^{n}\sum_{j=1}^{n}\left(  1-\left\vert u_{ij}(m_{\alpha})\right\vert
^{2}\right)  ^{1/2}\leq Cn\lambda\Vert A-A_{\alpha}\Vert_{\infty}.
\]
By taking limits on both sides of the above inequality, one obtains
\begin{align*}
\sum_{i=1}^{n}\sum_{j=1}^{n}\left(  1-\left\vert u_{ij}(m)\right\vert
^{2}\right)  ^{1/2}  & =\limsup_{m_{\alpha}\rightarrow m}\sum_{i=1}^{n}%
\sum_{j=1}^{n}\left(  1-\left\vert u_{ij}(m_{\alpha})\right\vert ^{2}\right)
^{1/2}\\
& \leq\limsup_{m_{\alpha}\rightarrow m}Cn\lambda\Vert A-A_{\alpha}%
\Vert_{\infty}=0.
\end{align*}
That is,
\[
\sum_{i=1}^{n}\sum_{j=1}^{n}\left(  1-\left\vert u_{ij}(m)\right\vert
^{2}\right)  ^{1/2}=0.
\]
On the other hand, it follows from (\ref{unimodular}) that
\begin{align*}
\limsup_{z\rightarrow m}\sum_{i=1}^{n}\left\Vert H_{\lambda\Phi(I-A)}\left(
k_{z}e_{i}\right)  \right\Vert  & \leq\limsup_{z\rightarrow m}\sum_{i=1}%
^{n}\sum_{j=1}^{n}\left(  1-\left\vert u_{ij}(z)\right\vert ^{2}\right)
^{1/2}\\
& =\sum_{i=1}^{n}\sum_{j=1}^{n}\left(  1-\left\vert u_{ij}(m)\right\vert
^{2}\right)  ^{1/2}=0.
\end{align*}
The proof is complete.
\end{proof}

\section{Some applications}

\label{Section 5}

To better illustrate Theorem \ref{main theorem TFAE}, we provide two examples
of compact products below. First, we present a new proof of the main results
from \cite{ChengChu15}, based on Theorem \ref{maina}. Second, we use
Corollary \ref{cor2} to emphasize the intricacies and difficulties associated
with Theorem \ref{maina}.

\begin{corollary}
[\cite{ChengChu15}]Let $\varphi,\psi\in L^{\infty}.$ Then $H_{\varphi}%
T_{\psi}$ is compact if and only if
\[
H^{\infty}\left[  \varphi\right]  \cap H^{\infty}\left[  \psi,\varphi
\psi\right]  \subset H^{\infty}+C.
\]

\end{corollary}

\begin{proof}
By Theorem \ref{maina}, $H_{\varphi}T_{\psi}$ is compact if and only
\[
\bigcap_{A\in\left(  M_{1}\right)  _{2^{2}}}H^{\infty}\left[  \varphi
(I-A),A\psi,\varphi A\psi\right]  \subset H^{\infty}+C.
\]

In this case, the matrix $A$ is a scalar. In the case $A=0,$ we get
$H^{\infty}\left[  \varphi(I-A),A\psi,\varphi A\psi\right]  =H^{\infty}\left[
\varphi\right]  $.

In the case $A=1,$ we get $H^{\infty}\left[  \varphi(I-A),A\psi,\varphi
A\psi\right]  =H^{\infty}\left[  \psi,\varphi\psi\right]  .$ In all other
cases, we get $H^{\infty}\left[  \varphi, \psi,\varphi\psi\right]  $. Thus we
will have
\[
H^{\infty}\left[  \varphi, \psi,\varphi\psi\right]  \cap H^{\infty}\left[
\varphi\right]  \cap H^{\infty}\left[  \psi,\varphi\psi\right]  \subset
H^{\infty}+C
\]
if and only if for each support set $s(m)$, one of the following holds: (1)
$\varphi|_{s(m)}\in H^{\infty}|_{s(m)}$, $\psi|_{s(m)}\in H^{\infty}|_{s(m)}$,
$(\varphi\psi)|_{s(m)}\in H^{\infty}|_{s(m)}.$ (2) $\varphi|_{s(m)}\in
H^{\infty}|_{s(m)}.$ (3) $\psi|_{s(m)}\in H^{\infty}|_{s(m)}$, $(\varphi
\psi)|_{s(m)}\in H^{\infty}|_{s(m)}.$ Therefore, we have
\[
H^{\infty}\left[  \varphi, \psi,\varphi\psi\right]  \cap H^{\infty}\left[
\varphi\right]  \cap H^{\infty}\left[  \psi,\varphi\psi\right]  \subset
H^{\infty}+C
\]
if and only if
\[
H^{\infty}\left[  \varphi\right]  \cap H^{\infty}\left[  \psi,\varphi
\psi\right]  \subset H^{\infty}+C,
\]
which completes the proof.
\end{proof}

\begin{corollary}
\label{cor2} Let $\varphi_{i},\psi_{i}\in L^{\infty}$ for $i=1,2.$ Then
$H_{\varphi_{1}}T_{\psi_{1}}+H_{\varphi_{2}}T_{\psi_{2}}$ is compact if and
only if
\begin{align}
&  H^{\infty}\left[  \varphi_{1},\varphi_{2}\right]  \cap H^{\infty}\left[
\psi_{1},\psi_{2},\varphi_{1}\psi_{1}+\varphi_{2}\psi_{2}\right]
\label{line1}\\
&  \cap H^{\infty}\left[  \varphi_{2},\psi_{1},\varphi_{1}\psi_{1}\right]
\label{line-1}\\
&  \cap H^{\infty}\left[  \varphi_{1},\psi_{2},\varphi_{2}\psi_{2}\right]
\label{line0}\\
&  \cap_{\lambda\neq0}H^{\infty}\left[  \varphi_{2}-\lambda\varphi_{1}%
,\psi_{1}+\lambda\psi_{2},\varphi_{1}\left(  \psi_{1}+\lambda\psi_{2}\right)
\right] \label{line2}\\
&  \cap_{\lambda\neq0}H^{\infty}\left[  \varphi_{2}-\lambda\varphi_{1}%
,\psi_{1}+\lambda\psi_{2},\varphi_{2}\left(  \psi_{1}+\lambda\psi_{2}\right)
\right] \label{line3}\\
&  \subset H^{\infty}+C.\nonumber
\end{align}

\end{corollary}

\begin{proof}
Write
\begin{align}
\Phi=\left[
\begin{array}
[c]{cc}%
\varphi_{1} & \varphi_{2}\\
0 & 0
\end{array}
\right]  ,\Psi=\left[
\begin{array}
[c]{cc}%
\psi_{1} & 0\\
\psi_{2} & 0
\end{array}
\right]  .\nonumber
\end{align}
By Theorem \ref{maina}, $H_{\varphi_{1}}T_{\psi_{1}}+H_{\varphi_{2}}%
T_{\psi_{2}}$ is compact if and only if
\[
\bigcap_{A\in\left(  M_{2}\right)  _{2^{4}}}H^{\infty}\left[  \Phi
(I-A),A\Psi,\Phi A\Psi\right]  \subset H^{\infty}_{2\times2}+C_{2\times2}.
\]

\textbf{{Case A:}} In the case $rankA=2$ or $rank(I-A)=2,$ we have%
\[
H^{\infty}\left[  \Phi(I-A),A\Psi,\Phi A\Psi\right]  = H^{\infty}\left[
\Psi,\Phi\Psi\right]  \text{ or }H^{\infty}\left[  \Phi(I-A),A\Psi,\Phi
A\Psi\right]  =H^{\infty}\left[  \Phi\right]  ,
\]
i.e.
\[
H^{\infty}\left[  \Phi(I-A),A\Psi,\Phi A\Psi\right]  =H^{\infty}\left[
\left[
\begin{array}
[c]{cc}%
\psi_{1} & 0\\
\psi_{2} & 0
\end{array}
\right]  , \left[
\begin{array}
[c]{cc}%
\varphi_{1}\psi_{1}+\varphi_{2}\psi_{2} & 0\\
0 & 0
\end{array}
\right]  \right]  \quad\text{or}\quad H^{\infty} \left[
\begin{array}
[c]{cc}%
\varphi_{1} & \varphi_{2}\\
0 & 0
\end{array}
\right]  ,
\]
which is the same Douglas algebras as (\ref{line1}).

\textbf{{Case B:}} Next, it suffices to discuss $A\in\left(  M_{2}\right)
_{16}$ such that $rankA=rank(I-A)=1.$ In this situation, there are three
cases:%
\begin{align}
A  &  =\left[
\begin{array}
[c]{cc}%
1 & \lambda\\
0 & 0
\end{array}
\right]  ,\label{line1a}\\
A  &  =\left[
\begin{array}
[c]{cc}%
0 & 0\\
\lambda & 1
\end{array}
\right]  ,\label{line2a}\\
A  &  =\left[
\begin{array}
[c]{cc}%
1-\lambda t & \lambda\left(  1-\lambda t\right) \\
t & \lambda t
\end{array}
\right]  ,\left\vert \lambda\right\vert \leq16. \label{line3a}%
\end{align}

\textbf{{Case 1:}} If $\lambda=0$ in the matrix (\ref{line1a}), we have
\[
H^{\infty}\left[  \Phi(I-A),A\Psi,\Phi A\Psi\right]  =H^{\infty} \left[
\left[
\begin{array}
[c]{cc}%
0 & \varphi_{2}\\
0 & 0
\end{array}
\right]  ,\left[
\begin{array}
[c]{cc}%
\psi_{1} & 0\\
0 & 0
\end{array}
\right]  ,\left[
\begin{array}
[c]{cc}%
\varphi_{1}\psi_{1} & 0\\
0 & 0
\end{array}
\right]  \right]  ,
\]
which is the same Douglas algebras as (\ref{line-1}).

If $\lambda\neq0$ in the matrix (\ref{line1a}),%
\begin{align}
\Phi(I-A)=\left[
\begin{array}
[c]{cc}%
\varphi_{1} & \varphi_{2}\\
0 & 0
\end{array}
\right]  \left[
\begin{array}
[c]{cc}%
0 & -\lambda\\
0 & 1
\end{array}
\right]  =\left[
\begin{array}
[c]{cc}%
0 & \varphi_{2}-\lambda\varphi_{1}\\
0 & 0
\end{array}
\right]  ,\nonumber
\end{align}
\begin{align}
A\Psi= \left[
\begin{array}
[c]{cc}%
1 & \lambda\\
0 & 0
\end{array}
\right]  \left[
\begin{array}
[c]{cc}%
\psi_{1} & 0\\
\psi_{2} & 0
\end{array}
\right]  =\left[
\begin{array}
[c]{cc}%
\psi_{1}+\lambda\psi_{2} & 0\\
0 & 0
\end{array}
\right]  ,\nonumber
\end{align}

\begin{align}
\Phi A\Psi=\left[
\begin{array}
[c]{cc}%
\varphi_{1} & \varphi_{2}\\
0 & 0
\end{array}
\right]  \left[
\begin{array}
[c]{cc}%
\psi_{1}+\lambda\psi_{2} & 0\\
0 & 0
\end{array}
\right]  =\left[
\begin{array}
[c]{cc}%
\varphi_{1}\left(  \psi_{1}+\lambda\psi_{2}\right)  & 0\\
0 & 0
\end{array}
\right]  .\nonumber
\end{align}
Therefore, we have
\[
H^{\infty}\left[  \Phi(I-A),A\Psi,\Phi A\Psi\right]  =H^{\infty} \left[
\left[
\begin{array}
[c]{cc}%
0 & \varphi_{2}-\lambda\varphi_{1}\\
0 & 0
\end{array}
\right]  ,\left[
\begin{array}
[c]{cc}%
\psi_{1}+\lambda\psi_{2} 0 & \\
0 & 0
\end{array}
\right]  ,\left[
\begin{array}
[c]{cc}%
\varphi_{1}\left(  \psi_{1}+\lambda\psi_{2}\right)  & 0\\
0 & 0
\end{array}
\right]  \right]  ,
\]
which is the same Douglas algebras as (\ref{line2}).

\textbf{{Case 2:}} If $\lambda=0$ in the matrix (\ref{line2a}), we have
\[
H^{\infty}\left[  \Phi(I-A),A\Psi,\Phi A\Psi\right]  =H^{\infty} \left[
\left[
\begin{array}
[c]{cc}%
\varphi_{1} & 0\\
0 & 0
\end{array}
\right]  ,\left[
\begin{array}
[c]{cc}%
0 & 0\\
\psi_{2} & 0
\end{array}
\right]  ,\left[
\begin{array}
[c]{cc}%
\varphi_{2}\psi_{2} & 0\\
0 & 0
\end{array}
\right]  \right]  ,
\]
which is the same Douglas algebras as (\ref{line0}).

If $\lambda\neq0$ in the matrix (\ref{line2a}),
\begin{align}
\Phi(I-A)=\left[
\begin{array}
[c]{cc}%
\varphi_{1} & \varphi_{2}\\
0 & 0
\end{array}
\right]  \left[
\begin{array}
[c]{cc}%
1 & 0\\
-\lambda\  & 0
\end{array}
\right]  =\left[
\begin{array}
[c]{cc}%
-\lambda\left(  \varphi_{2}-\frac{1}{\lambda}\varphi_{1}\right)  & 0\\
0 & 0
\end{array}
\right]  ,\nonumber
\end{align}

\begin{align}
A\Psi= \left[
\begin{array}
[c]{cc}%
0 & 0\\
\lambda & 1
\end{array}
\right]  \left[
\begin{array}
[c]{cc}%
\psi_{1} & 0\\
\psi_{2} & 0
\end{array}
\right]  =\left[
\begin{array}
[c]{cc}%
0 & 0\\
\lambda\left(  \psi_{1}+\frac{1}{\lambda}\psi_{2}\right)  & 0
\end{array}
\right]  ,\nonumber
\end{align}

\begin{align}
\Phi A\Psi &  =\left[
\begin{array}
[c]{cc}%
\varphi_{1} & \varphi_{2}\\
0 & 0
\end{array}
\right]  \left[
\begin{array}
[c]{cc}%
0 & 0\\
\lambda\left(  \psi_{1}+\frac{1}{\lambda}\psi_{2}\right)  & 0
\end{array}
\right]  =\left[
\begin{array}
[c]{cc}%
\lambda\varphi_{2} \left(  \psi_{1}+\frac{1}{\lambda}\psi_{2}\right)  & 0\\
0 & 0
\end{array}
\right]  .\nonumber
\end{align}
Therefore, we have \begin{flalign}
&H^{\infty}\left[  \Phi(I-A),A\Psi,\Phi A\Psi\right] \nonumber
\\ \nonumber&=H^{\infty} \left[ \left[
\begin{array}
[c]{cc}%
-\lambda\left( \varphi_{2}-\frac{1}{\lambda}\varphi_{1}\right)  & 0  \\
0 & 0
\end{array}
\right] ,\left[
\begin{array}
[c]{cc}%
0 & 0  \\
\lambda\left( \psi_{1}+\frac{1}{\lambda}\psi_{2}\right)  & 0
\end{array}
\right] ,\left[
\begin{array}
[c]{cc}%
\lambda\varphi_{2} \left( \psi_{1}+\frac{1}{\lambda}\psi_{2}\right)  & 0  \\
0 & 0
\end{array}
\right] \right] ,
\end{flalign}
which is the same Douglas algebras as (\ref{line3}).

\textbf{{Case 3:}} We will show that the matrix (\ref{line3a}) does not yield
new Douglas algebras.

\textbf{{Case 3.1:}} When $t=0$, the matrix $A$ is (\ref{line1a}). This case
is reduced to \textbf{{Case 1}}

\textbf{{Case 3.2:}} When $1-\lambda t=0$, the matrix $A$ is $\left[
\begin{array}
[c]{cc}%
0 & 0\\
\frac{1}{\lambda} & 1
\end{array}
\right]  $, and we get the same Douglas algebras as the case (\ref{line2a}).
This case is reduced to \textbf{{Case 2}}

\textbf{{Case 3.3:}} When $t\neq0$ and $1-\lambda t\neq0$, the matrix $A$
yields the linear combination of (\ref{line2}) and (\ref{line3}), due to the
following equations.%

\begin{align}
\Phi(I-A)=\left[
\begin{array}
[c]{cc}%
\varphi_{1} & \varphi_{2}\\
0 & 0
\end{array}
\right]  \left[
\begin{array}
[c]{cc}%
\lambda t & -\lambda\left(  1-\lambda t\right) \\
-t & 1-\lambda t
\end{array}
\right]  =\left[
\begin{array}
[c]{cc}%
-t\left(  \varphi_{2}-\lambda\varphi_{1}\right)  & \left(  1-\lambda t\right)
\left(  \varphi_{2}-\lambda\varphi_{1}\right) \\
0 & 0
\end{array}
\right]  ,\nonumber
\end{align}

\begin{align}
A\Psi= \left[
\begin{array}
[c]{cc}%
1-\lambda t & \lambda\left(  1-\lambda t\right) \\
t & \lambda t
\end{array}
\right]  \left[
\begin{array}
[c]{cc}%
\psi_{1} & 0\\
\psi_{2} & 0
\end{array}
\right]  =\left[
\begin{array}
[c]{cc}%
\left(  1-\lambda t\right)  \left(  \psi_{1}+\lambda\psi_{2}\right)  & 0\\
t\left(  \psi_{1}+\lambda\psi_{2}\right)  & 0
\end{array}
\right]  ,\nonumber
\end{align}

\begin{align}
\Phi A\Psi=\left[
\begin{array}
[c]{cc}%
\varphi_{1} & \varphi_{2}\\
0 & 0
\end{array}
\right]  \left[
\begin{array}
[c]{cc}%
\left(  1-\lambda t\right)  \left(  \psi_{1}+\lambda\psi_{2}\right)  & 0\\
t\left(  \psi_{1}+\lambda\psi_{2}\right)  & 0
\end{array}
\right]  =\left[
\begin{array}
[c]{cc}%
\left(  1-\lambda t\right)  \varphi_{1}\left(  \psi_{1}+\lambda\psi
_{2}\right)  +t\varphi_{2}\left(  \psi_{1}+\lambda\psi_{2}\right)  & 0\\
0 & 0
\end{array}
\right]  .\nonumber
\end{align}
Combing all the above results, we complete the proof of the Corollary.
\end{proof}

\noindent\textbf{Acknowledgement}: P. Ma was supported by NNSF of China (Grant
numbers 12171484,12571140), the Natural Science Foundation of Hunan
Province (Grant number 2023JJ20056), the Science and Technology Innovation
Program of Hunan Province (Grant number 2023RC3028), and Central South
University Innovation-Driven Research Programme (Grant number 2023CXQD032). M.
Li was also supported by the Graduate Research and Innovation Projects of Hunan Province (Grant number CX20250159).

Caixing Gu

Department of Mathematics, California Polytechnic State University, San Luis
Obispo, CA 93407, USA

E-mail: cgu@calpoly.edu\bigskip\bigskip

Meng Li

School of Mathematics and Statistics, Hunan Research Center of the Basic Discipline for Analytical Mathematics, HNP-LAMA, Central South University, Changsha 410083, China

E-mail address: mengli97@csu.edu.cn\bigskip\bigskip

Pan Ma

School of Mathematics and Statistics, Hunan Research Center of the Basic Discipline for Analytical Mathematics, HNP-LAMA, Central South University, 
Changsha 410083, China

E-mail address: pan.ma@csu.edu.cn

\bigskip

Mathematics Subject Classification (2020). 47B35, 46J15, 46J20, 30H10

\bigskip

\textbf{Keywords:} Hardy space, block Hankel operators, block Toeplitz
operators, compact operators
\end{document}